\documentclass[12pt]{amsart}

\usepackage{graphicx,lipsum,pgfplots}
\usepgfplotslibrary{fillbetween}
\usepackage{float}
\pgfplotsset{width=10.6cm,compat=1.16}

\usepackage{geometry,a4wide,amssymb,amsfonts,amsrefs,amsmath}         
\usepackage{bbm,xcolor}
\usepackage{dsfont}
\usepackage{yfonts}
\usepackage[T1]{fontenc}
\usepackage{mathtools}
\usepackage{mathrsfs}

\newcommand{\miniscule}{\@setfontsize\miniscule{8}{5}}

\newcommand\mystyle{\everymath{\displaystyle}}
\mystyle

\newcommand{\Z}{\mathbb{Z}}
\newcommand{\A}{\mathcal{A}}
\newcommand{\N}{\mathbb{N}}
\newcommand{\R}{\mathbb{R}}
\newcommand{\F}{\mathcal{F}}
\newcommand{\G}{\mathscr{G}}

\newcommand{\FP}{\mathfrak{P}}
\newcommand{\La}{\mathcal{L}}
\newcommand{\bL}{\mathbb{L}}

\newcommand{\cov}{\text{cov}}
\newcommand{\gl}{\overset{\rightharpoonup}{\gamma}}
\newcommand{\gPhi}{\overset{\rightharpoonup}{\Phi}}
\newcommand{\glr}{\overset{\rightleftharpoons}{\gamma}}

\def\mybf #1{\textbf{\textit{#1}}}

\theoremstyle{plain}

\newtheorem{thm}{Theorem}[section]
\newtheorem{lem}[thm]{Lemma}
\newtheorem{prop}[thm]{Proposition}

\newtheorem{defin}[thm]{Definition}

\newtheorem{example}[thm]{Example}
\newtheorem{remark}[thm]{Remark}

\newtheorem{manualtheoreminner}{Theorem}
\newenvironment{manualtheorem}[1]{%
  \renewcommand\themanualtheoreminner{#1}%
  \manualtheoreminner
}{\endmanualtheoreminner}

\usepackage{setspace}
\usepackage{fancyhdr}

\usepackage{ulem}

\title[On an extension of a theorem by Ruelle to long-range potentials]{On an extension of a theorem by Ruelle to long-range potentials}

\author[A.C.D. van Enter]{Aernout C.D. van Enter}
\address[Aernout van Enter]{Johann Bernoulli Institute for Mathematics and Computer Science, 
University of Groningen, Nijenborgh 9, 9747 AG Groningen, Netherlands}
\email{a.c.d.van.enter@rug.nl}
\author[R. Fern\'andez]{Roberto Fern\'andez}
\address[Roberto Fern\'andez]{New York University Shanghai, 1555 Century Avenue, Pudong, Shanghai, China.}
\email{rf87@nyu.edu}

\author[M. Makhmudov]{Mirmukhsin Makhmudov}
\address[Mirmukhsin Makhmudov]{Mathematical Institute,  Leiden University, Einsteinweg 55, 2333 CC Leiden, The Netherlands}
\email{m.makhmudov@math.leidenuniv.nl}
\author[E. Verbitskiy]{Evgeny Verbitskiy}
\address[Evgeny Verbitskiy]{Mathematical Institute, Leiden University, Einsteinweg 55, 2333 CC Leiden, The Netherlands, and} 
\address{Korteweg-de Vries Institute for Mathematics, University of Amsterdam, 
 Science Park 105-107, 1098 XG Amsterdam, The Netherlands
}
\email{evgeny@math.leidenuniv.nl}

\begin{document}

\begin{abstract} Ruelle's transfer operator plays an important role in understanding thermodynamic and probabilistic properties of dynamical systems. In this work, we develop a method of finding eigenfunctions of transfer operators based on comparing Gibbs measures on the half-line $\mathbb Z_+$  and the whole line $\Z$.
For a rather broad class of potentials, including both the ferromagnetic and antiferromagnetic long-range Dyson potentials, we are able to establish the existence of integrable, but not necessarily continuous, eigenfunctions.
For a subset thereof we prove that the eigenfunction is actually continuous.
\end{abstract}
\dedicatory{To Geoffrey Grimmett, with appreciation for all he has provided to us, both personally and probabilistically,  by his presence among the percolation and correlated percolation practitioners.}
\maketitle
\section{Introduction}
 One of the main issues of equilibrium statistical mechanics is to derive and describe the  properties of the possible (global) states of a macroscopic system
starting from the
knowledge of the finite-volume (local) states of the system.
In order to give a mathematical framework to this problem, R. Dobrushin, O. Lanford, and D. Ruelle developed the so-called DLR formalism in the second half of the last century. 

Shortly after their introduction of Gibbs measures,
Spitzer and Averintsev \cites{Av,Spi}  characterised Gibbs measures for short-range potentials in terms of measures having Markov properties. 
These results then were simplified and generalised, in various directions by Hammersley and Clifford {\cite{HC1968}}, by Sullivan, Kozlov {\cites{Kozlov1974,Sullivan1973}}, and by Geoffrey Grimmett in his first paper \cite{Gri1}.

The novel DLR formalism was immediately adapted to the theory of  Dynamical Systems by Ya. Sinai \cites{Sin,Sinai-Uspekhi,Sinai-ICM}.
There are, however, two important differences between the models typically studied in Statistical Mechanics and those studied in Dynamical Systems. 
First, in Dynamical Systems one is typically interested in one-dimensional systems, that is, systems
with configuration spaces $E^{\Z}$, where $E$ is the set of possible spin values and the spatial dimension 
represents time. The conditional probabilities in Dynamical Systems are typically "one-sided" (from past to future), those in Statistical Mechanics "two-sided" (from outside to inside). For short-range potentials this does not make much of a difference, although in general regularity properties between one-sided and two-sided conditional probabilities may differ \cites{BFV2019,  BEvELN,ELP, ES,FGM}.
Second, and perhaps more important, the natural description of dynamical systems often  involves 
{\it half-line}
configuration spaces $E^{\Z_+}$, rather than { \it whole-line} 
configuration spaces of the form $E^{\Z}$. 

Already in his original papers, Sinai addressed these questions, \cites{Sinai-ICM,Sinai-Uspekhi,Sin}. 
He showed that 
  Gibbs equilibrium states for \emph{exponentially decaying interactions} 
are, in fact, equilibrium states for half-line potentials as well.
The issue of half-line versus whole-line Gibbsianness is, therefore, as old as the theory of thermodynamic formalism in Dynamical Systems.
For more recent results established in this area see, e.g., \cites{Walters1975, Walters1978, Walters2001, Walters2005, FM2004,
BFV2019}.
One of the most important Dynamical Systems tools used in the study of Gibbs equilibrium states is the
so-called Ruelle's transfer operator \cites{Rue, Walters2001}. Various probabilistic properties
of chaotic dynamical systems can be characterized in terms of transfer
 operators. In particular, a central issue is to characterize those potentials 
 (interactions) for which transfer operators have positive continuous 
 eigenfunctions that is,  finding those potentials for which Ruelle's theorem holds.
This question has been answered by different authors \cites{Walters2001, Walters1978, Walters1975, Fan1995, Ruelle1968, Mayer} for different regularity classes of potentials. 
In this paper, we answer it for potentials beyond these earlier studied classes.

Our paper is organized as follows: 
\begin{itemize}
    \item In Section \ref{section about specifications},  Section \ref{Ruelle's theorem  and Equilibrium states},  and Section \ref{Basic notions III:  Equivalence of interactions.  Dobrushin uniqueness condition} we introduce the notions of thermodynamic formalism that are important for this paper.
    \item In Section \ref{novel class of potentials}, 
    we discuss the relationship between half-line and whole-line Gibbs measures and formulate the first part of the main results (Theorem \ref{convergence of shifted measures} and 
    \ref{convergence of averaged measures}).
    \item In Section \ref{Eigenfunctions of transfer operators and absolute continuity of Gibbs measures}, we discuss when whole-line Gibbs measures are absolutely continuous with respect to the product  of two half-line ones, and we formulate the second series of our main results (Theorem \ref{main theorem about the Dyson model}, 
    \ref{UI and Dunford-Pettis} 
    and
    \ref{sufficient cond for UI}).
    \item In Section \ref{Applications: Dyson model}, the Dyson model, the main example of the paper, is discussed. 
    We then discuss what is the behaviour in other regimes of the phase transitions.
    \item Section \ref{proofs of the main results} and Section  \ref{final remarks} are dedicated to the proofs of our main results and the final remarks. 
\end{itemize}

\section{Basic notions: I. Specifications}\label{section about specifications}

The  Dobrushin-Lanford-Ruelle definition of Gibbs states via specifications goes well beyond the standard lattices $\Z^d$.
Consider the \emph{lattice system} $\Omega=E^{\mathbb L}$, where $\mathbb L$ is an at most countable set (lattice)
and $E$ is a set of possible spin values. 
In this paper, we will focus on \emph{finite} $E$.
We denote the Borel $\sigma-$algebra of the measurable subsets of $\Omega$ by $\F$.
The specification is a \emph{consistent} family of probability kernels (conditional probabilities) indexed by finite subsets $\Lambda$ of $\mathbb L$ denoted by $\Lambda\Subset \mathbb L$. 
The consistency condition is the requirement that 
$\gamma_{\Lambda}\gamma_{\Delta}=\gamma_{\Lambda}$ for all $\Delta\subset\Lambda\Subset \mathbb L$ 
\cite[Chapter 1]{Georgii-book} .

 In the sequel, for $\Delta,\Lambda\subset \bL$,   we will denote 
 the concatenation of strings $\xi_{\Delta}\in E^{\Delta}$, $\eta_{\Lambda}\in E^{\Lambda}$ by $\xi_{\Delta}\eta_{\Lambda}$, namely, $\xi_{\Delta}\eta_{\Lambda}$ is a string such that $(\xi_{\Delta}\eta_{\Lambda})_{i}=\xi_i$ if $i\in \Delta$ and $(\xi_{\Delta}\eta_{\Lambda})_{i}=\eta_i$ if $i\in \Lambda$. 
 Given the specification 
$\gamma=\{\gamma_\Lambda\}_{\Lambda\Subset \mathbb L}$ on $\Omega=E^\bL$, we say that a probability measure $\mu$ on $\Omega$ is Gibbs for the specification 
$\gamma$ or, equivalently,  that $\mu$ is consistent with $\gamma$, if 
$$
\mu(\sigma_\Lambda|\sigma_{\Lambda^c}) = \gamma_{\Lambda}(\sigma_\Lambda\sigma_{\Lambda^c})\quad 
\text{for $\mu$-a.a. }\sigma\in \Omega,
$$
or, equivalently, if the DLR equations hold:
$$
\int \gamma_\Lambda f \, d\mu = \int f\, d\mu,
$$
for all $f\in L^1(\Omega,\mu)$, and every $\Lambda\Subset \mathbb L$, where
$$
\gamma_\Lambda f(\sigma) = \sum_{\xi_\Lambda\in E^\Lambda} \gamma_\Lambda(\xi_\Lambda|\sigma_{\Lambda^c})
f(\xi_\Lambda\sigma_{\Lambda^c}),
$$
note that for $\omega\in\Omega$ and $\Lambda\subset\bL$, $\omega_{\Lambda^c}$ denotes the (infinite) string $\omega_{\bL\setminus\Lambda}$.
The set of all Gibbs measures for $\gamma$ will be denoted by $\G(\Omega, \gamma)$.
For any probability measure $\mu$ one can find at least one specification $\gamma$ such that
$\mu$ is Gibbs for $\gamma$ \cite{Goldstein1978}.
However, useful and interesting specifications have additional properties
such as finite energy (non-nullness) and quasi-locality (continuity). 
We now turn to two particular ways of defining specifications as used in Statistical Mechanics and Dynamical Systems.

\subsection{Gibbs(ian) specifications in Statistical Mechanics}

 An interaction  $\Phi$ is a family of  functions $\{\Phi_{\Lambda}\}$,  indexed by finite subsets 
 $\Lambda\Subset \bL$, such that each function $\Phi_{\Lambda}$,
  depends only on values of $\sigma$ in $\Lambda$, that is,  with a slight abuse of notation, $\Phi_{\Lambda}(\sigma)=\Phi_\Lambda(\sigma_\Lambda)$. One needs to impose some additional summability conditions on 
 the interaction $\Phi$: $\Phi$ is said to be \mybf{uniformly absolutely convergent} (UAC) if  $\sup_{i\in\bL} \sum_{V\Subset \bL, V\ni i}||\Phi_{V}||_{\infty}< \infty$. 
For an UAC interaction $\Phi$, the specification (specification density)   $\gamma^{\Phi}=\{\gamma^{\Phi}_\Lambda\}_{\Lambda\Subset\bL}$ is defined as follows, for  $\omega,\eta\in\Omega$,
\begin{equation}\label{spec density for UAC interaction}
    \gamma^{\Phi}_{\Lambda}(\omega_\Lambda|\eta_{\Lambda^c}):=\frac{e^{-H_{\Lambda}^{\Phi}(\omega_{\Lambda}\eta_{\Lambda^c})}}{Z^{\Phi}_{\Lambda}(\eta)},
\end{equation}
where 
$H_{\Lambda}^{\Phi}(\omega):=\sum_{V\cap \Lambda\neq \emptyset}\Phi_V(\omega)$ is the Hamiltonian in the volume $\Lambda$, and $Z_{\Lambda}^{\Phi}$ is a normalization constant (the partition function), i.e., $Z_{\Lambda}^{\Phi}(\eta):=\sum_{\bar \omega_{\Lambda}\in E^{\Lambda}}e^{-H_{\Lambda}^{\Phi}(\bar \omega_{\Lambda}\eta_{\Lambda^c})}$.

It should be stressed that a Gibbsian specification $\gamma^{\Phi}$ is always \textit{quasilocal} \cites{BGMMT2020,VEFS1993,Georgii-book}. 
In the current setting, in which $E$ is finite, this property is equivalent to the fact that for all $\Lambda\Subset \bL$ and $\omega_{\Lambda}\in E^{\Lambda}$, $\gamma_{\Lambda}(\omega_{\Lambda}|\eta)$ is a continuous function of the boundary condition $\eta\in\Omega$. 
Another important property of the Gibbsian specifications is   \textit{non-nullness}, which means that for all  volumes $\Lambda\Subset \bL$, 
$\inf_{\eta, \omega\in\Omega}\gamma^{\Phi}_{\Lambda}(\omega_{\Lambda}|\eta_{\Lambda^c})>0$.

We denote the set of Gibbs states for the interaction $\Phi$ by $\G(\Omega, \Phi)$ (or $\G(\Phi)$).  It is a convex set ---in fact a simplex--- which is always non-empty if,  as is the case in this paper, the spin space $E$ is compact.

Depending on the symmetries of the lattice $\bL$ and the spin space $E$, the interactions and specifications may also exhibit some symmetries. For example, if $\bL=\Z$,  then an interaction  $\Phi$ on $X:=E^\Z$ is called \textbf{\textit{translation-invariant}} if for all $\Lambda\Subset\Z$,
every $k\in\Z$ and $\omega\in X$, $\Phi_{\Lambda+k}(\omega)=\Phi_\Lambda(S^k(\omega))$, where $\Lambda+k:=\{i+k:i\in\Lambda\}$.
Respectively, a specification $\gamma$ on $X=E^\Z$ is called \textbf{\textit{translation-invariant}} if for all $B\in \F$, $\Lambda\Subset\Z$, $k\in\Z$ and $\omega\in X$, $\gamma_{\Lambda+k}(B|\omega)=\gamma_{\Lambda}(S^k(B)|S^k(\omega))$. 
Translation-invariant interactions give rise to translation-invariant specifications. 
To some extent, the opposite statement is also true: Sullivan showed (\cites{Sullivan1973, BGMMT2020}) that 
for a quasilocal translation-invariant specification on $\Z$, one can find a translation-invariant 
interaction $\Phi$ such that $\gamma=\gamma^{\Phi}$. Recently, however, it was shown in \cite{BGMMT2020} that  this interaction
is not necessarily uniformly absolutely convergent. 

It should be noted that Gibbsian specifications can be uniquely recovered from a consistent family of
single-site probability kernels (densities) $\{\gamma_{\{i\}}:\, i\in\mathbb L\}$ \cite{FM2004}. 
Therefore, it is sufficient to study only the single-site densities of a Gibbsian specification instead of studying {\it all} densities. 
Due to this fact, it is worth defining the single-site densities of a specification separately from the concept of specification as follows.

\begin{defin}\label{def single site densit}\cite{F2006}
A collection $\{\gamma_{\{i\}}\}_{i\in \bL}$ of positive functions $\gamma_{\{i\}}(\cdot|\cdot):E\times E^{\bL\setminus \{i\}}\to (0,1)$ is called the \textit{family of single-site densities of a specification} if
\begin{itemize}
    \item[(i)]$\sum_{a_i\in E}\gamma_{\{i\}}(a_i|\omega_{\{i\}^c})=1$ for all $\omega\in\Omega=E^\bL$, $i\in \bL$,
    \item[(ii)] and for all $i,j\in \bL$, $\alpha,\omega\in\Omega$ the following holds $$\frac{\gamma_{\{i\}}(\alpha_i|\alpha_j\omega_{\{i,j\}^c})}{\sum_{\beta_{\{i,j\}}}\frac{\gamma_{\{j\}}(\beta_j|\beta_i\omega_{\{i,j\}^c})\gamma_{\{i\}}(\beta_i|\alpha_j\omega_{\{i,j\}^c})}{\gamma_{\{j\}}(\alpha_j|\beta_i\omega_{\{i,j\}^c})}}=
\frac{\gamma_{\{j\}}(\alpha_j|\alpha_i\omega_{\{i,j\}^c})}{\sum_{\beta_{\{i,j\}}}\frac{\gamma_{\{i\}}(\beta_i|\beta_j\omega_{\{i,j\}^c})\gamma_{\{j\}}(\beta_j|\alpha_i\omega_{\{i,j\}^c})}{\gamma_{\{i\}}(\alpha_i|\beta_j\omega_{\{i,j\}^c})}}.$$
\end{itemize}
\end{defin}
The following theorem signifies the importance of single-site densities of a specification, and it will be useful later. 
\begin{prop}\label{important theorem about single-site densities}\cite{FM2004}
Let $\{\gamma_{\{i\}}\}_{i\in \bL}$ be single-site densities of a specification. 
There is a unique non-null specification $\gamma$ on $(\Omega,\mathcal{F})$ having $\{\gamma_{\{i\}}\}_{i\in \bL}$ as its single-site densities. 
Furthermore, $\gamma$ is quasilocal if and only if all functions in the collection $\{\gamma_{\{i\}}\}_{i\in \bL}$ are continuous, and a probability measure $\mu\in\mathcal{M}_1(\Omega,\mathcal{F})$ is consistent with $\gamma$ if and only if it is consistent with all single-site probability kernels $\gamma_{\{i\}}$, ${i\in \bL}$.
\end{prop}


\subsection{Gibbs(ian) specifications in Dynamical Systems and Transfer operators}\label{Dynamical specifications}
As already mentioned above, in Dynamical Systems, the '\emph{natural}' lattice is the half-line  $\mathbb L=\mathbb Z_+$.
Let $X_+=E^{\Z_+}$ be the space of one-sided sequences $\omega=(\omega_n)_{n\ge 0}$ in alphabet $E$.
We equip $X_+$ with the metric $d(\omega,\omega') = \sum_{n=0}^\infty \mathbb I[\omega_n\ne\omega'_n] 2^{-n}$.
We also define the left-shift $S$ on $X_+$ by $y=Sx$ where $y_i=x_{i+1}$ for all $i\geq 0$.
Let $\phi:X_+\to\R$ be a continuous function (potential). Following \cites{CL,CLS,Walters2001},  we define the corresponding specification $\gl^{\phi}:=\{\gl_n=\gl_{[0,n-1]}^{\phi},\ n\ge 1\},$ 
by
\begin{equation}\label{the half-line specifications}
    \gl_n (a_0^{n-1}|x_{n}^\infty)= \frac {\exp\left( (S_n \phi)(a_0^{n-1}x_n^\infty)\right)}{\sum_{\bar a_0^{n-1}}
\exp\left( (S_n \phi)(\bar a_0^{n-1}x_n^\infty)\right)},\quad\text{where \/}
(S_n \phi)(x)=\sum_{k=0}^{n-1}\phi( S^k x).
\end{equation}
This gives a family of probability kernels on finite intervals $[0,n-1]$ in $\Z_+$. However, the definition extends to general  volumes $\Lambda\Subset \Z_+$ by 
\begin{equation}\label{the extension of the half-line specifications}
   \gl_{\Lambda}(a_{\Lambda}|x_{\Lambda^c})=\frac{\exp{((S_{n+1}\phi)(a_{\Lambda}x_{\Lambda^c}))}}{\sum_{\bar a_{\Lambda}}
\exp{((S_{n+1}\phi)(\bar a_{\Lambda}x_{\Lambda^c}))}}, \;\; a_{\Lambda}\in E^{\Lambda}, x\in X_+, 
\end{equation}
where $n=\max\Lambda$.

For $f\in C(X_+,\R)$,  one has
\begin{equation}\label{onesidedkernel}
\gl_n (f)(x) = \frac {\sum_{\bar a_0^{n-1}} \exp\left( (S_n \phi)(a_0^{n-1}x_n^\infty)\right) f(\bar a_0^{n-1} x_{n}^{\infty}) }{\sum_{\bar a_0^{n-1}}
\exp\left( (S_n \phi)(\bar a_0^{n-1}x_n^\infty)\right)}.
\end{equation}
It turns out that $\gl_n (f)$ can naturally be expressed in terms of the Ruelle-Perron-Frobenius transfer operator. This is the operator $\mathcal L_\phi$ acting on the space of continuous functions $C(X_+,\R)$ as
\begin{equation}\label{transfer}
\mathcal L_\phi f(x) =\sum_{y\in S^{-1}x} e^{\phi(y)} f(y)=\sum_{a\in E} e^{\phi(ax)} f(ax),
\end{equation}
where the configuration $ax$ is obtained by the concatenation of the letter $a$ and the configuration $x$.
Thus, for any  $n\ge 1$, 
$$
\mathcal L^n_\phi f(x)=\sum_{a_0^{n-1}\in E^n} e^{S_n\phi(a_0^{n-1}x)} f(a_0^{n-1}x),
\text{ and hence, }
\gl_n (f)(x) = \frac{ \mathcal L^n_\phi f(S^n x)} { \mathcal L^n_\phi \mathbf 1(S^n x)}.
$$
One readily checks that the family of probability kernels
$\{\gl_n\}$ has the standard properties of specifications;
most importantly, the consistency condition 
$$
\gl_m( \gl_n(f)) = \gl_n( \gl_m(f)) = \gl_m(f)
$$
for all $f\in C(X_+,\R)$  and every $m\ge n\ge 1$, and this particular property can be readily validated using properties of transfer operators \cite[Theorem 2.1]{Walters2001}.

The set of all Gibbs states on $X_+$ for potential $\phi$, i.e., the set of measures consistent with the specification 
$\gl^\phi$, will be denoted by $\G(X_+,\phi)$. 
The set of Gibbs measures $\G(X_+,\phi)$ is a closed convex set and the extremal points
of $\G(X_+,\phi)$ are tail-trivial.

Transfer operators allow for a \emph{dual view} on  Gibbs measures $\G(X_+,\phi)$. Define the dual  operator  $\mathcal L_\phi^*$,  acting on the space 
of  measures $\mathcal M(X_+)$ by

$$\int f\, d (\mathcal L_\phi^*\nu) = \int\mathcal L_\phi f\, d\nu\quad \text{ for all }f\in C(X_+,\R).
$$
It is well-known that there exists at least one eigenprobability $\nu$ on $X_+$ 
for the maximal eigenvalue $\lambda=e^{P(\phi)}$, i.e., 
$$
\mathcal L_\phi^*\nu=e^{P(\phi)} \nu, 
$$
where $P(\phi)=\lim_{n\to\infty} \frac 1n \log \sup_{x\in X_+} \mathcal L^n\mathbf 1(x)$ is the 
so-called topological pressure of $\phi$.
Combining the results of  \cite[Corollary 2.3]{Walters2001}  and \cite[Theorem 4.8]{CLS}, we can conclude that 
the sets of probability eigenmeasures of $\mathcal L_\phi^*$ and the Gibbs states for $\phi$ on $X_+$ coincide:
$$
\nu \in\G(X_+,\phi) \text{ if and only if \/} \mathcal L_\phi^*\nu=\lambda \nu.
$$

There is an interesting phenomenon in the theory of Gibbs measures on one-sided (half-line) symbolic spaces, which has no direct analogue in the two-sided (whole-line) context.  
Let us say that a continuous potential $\phi:X_+\to\R$ is quasi-normalized if $\mathcal{L}_{\phi} \mathbf 1$ is a constant function on $X_+$. 
It turns out that $\phi$ is quasi-normalized if and only if all Gibbs measures in $\G(X_+,\phi)$ are translation invariant.

\subsection{Relation between Gibbsian specifications}

Specifications discussed above are defined on different spaces: $X=E^{\Z}$ vs $X_+=E^{\Z_+}$,
as well as, in different terms: namely, the interaction $\Phi$ vs the potential $\phi$.
What is the relation between these classes of specifications?

For a given interaction $\Phi$, the potential $\phi$ should
be interpreted as minus the contribution to the energy from (the 
neighborhood of) the origin \cite[Section 3.2]{Rue}, \cite[Section 2.4.5]{VEFS1993}.
In fact there are multiple possibilities to define relevant $\phi$,
e.g.,
\begin{equation}\label{phi to Phi associationfull}
    \phi(\omega):=-\sum_{0\in V\Subset\Z}\frac {1}{|V|}\Phi_V(\omega_V),
\end{equation}
or,
\begin{equation}\label{phi to Phi association}
    \phi(\omega):=-\sum_{0\in V\Subset\Z_+}\Phi_V(\omega_V).
\end{equation}

In the setup of Thermodynamic Formalism, the second choice is particularly convenient, and it will be used in this paper. 

It is worth mentioning that for every $\phi\in C(E^{\Z_+},\R)$ there exists a translation-invariant interaction $\Phi^{\phi}$ on $X=E^{\Z}$ satisfying (\ref{phi to Phi association}) such that 
$$
\sum_{0\in V\Subset \Z}\frac{1}{|V|}||\Phi^{\phi}_V||_{\infty}<\infty,
$$
nonetheless, such $\Phi^{\phi}$ does not need to be unique \cite{Rue}.
The reciprocal property
 would be more interesting, namely that for each $\phi\in C(X_+)$ there would  exist
 a translation-invariant UAC interaction $\Phi$ on $X$ satisfying (\ref{phi to Phi association}).  
 Unfortunately, there exist counterexamples showing this to be false (c.f. Proposition \ref{phi-to-Phi_the form of corresponding whole-line spec} and Section \ref{final remarks}).

\section{Basic notions II:  Ruelle's theorem  and Equilibrium states}\label{Ruelle's theorem  and Equilibrium states}

Gibbs measures are eigenmeasures of the duals of transfer operators corresponding to the maximal eigenvalue. 
What can we say about the eigenfunctions of transfer operators? 
Ruelle established the first result
for \emph{smooth potentials} \cites{Ruelle1968, Rue}: if $\phi:X_+\to\R$ is H\"older continuous, then the transfer
operator $\mathcal L_\phi$ has a positive continuous eigenfunction $h$; $\mathcal L_\phi h=\lambda h$, with $\lambda=e^{P(\phi)}$. The existence of continuous eigenfunctions has been established for larger classes
of smooth potentials: for example, by Walters for potentials with summable variations  (\cites{Walters1978, Walters1975}). For less regular potentials -- the so-called Bowen class -- Walters
established the existence of bounded measurable eigenfunctions.

There is an interesting \emph{principal} relation between the Gibbs measures $\G(X_+,\phi)$, eigenfunctions of transfer operators,
and translation invariant equilibrium states.

\begin{prop}\label{existence of eigenfunction-criteria}
  Consider a continuous potential $\phi\in C(X_+)$, and suppose 
   $\nu\in\G(X_+,\phi)$ or, equivalently, $\mathcal L_\phi^* \nu =e^{P(\phi)}\nu$.
Then
   \begin{itemize}
       \item[1)] If there exists a non-negative eigenfunction $h\in L^1(X_+,\nu)$ of the transfer operator $\La_{\phi}$ with
       $\mathcal L_\phi h =e^{P(\phi)}h$, 
        then $\mu=h\cdot\nu$ (i.e., $d\mu=hd\nu$) is a translation invariant equilibrium state of $\phi$:
        namely, $\mu$ is the measure satisfying the variational principle
        $$
          h(S,\mu)+\int \phi d\mu =\sup_{\rho\in \mathcal M_1(X_+,S)} \Bigl[ h(S,\rho)+\int \phi d\rho\Bigr]=:P(\phi).
        $$
        where  $h(S,\cdot)$ is the Kolmogorov-Sinai entropy, the supremum is taken over the set of all 
        translation invariant probability measures on $X_+$, and $P(\phi)$ is the (topological) pressure of $\phi$.
       \item[2)] If there exists a translation invariant measure $\mu\in\mathcal{M}_1(X_+,S)$ such that $\mu\ll \nu$, then  $\mu$ is an equilibrium state for $\phi$  and the Radon-Nikodym derivative $h=\frac{d\mu}{d\nu}$ is the eigenfunction of $\La_{\phi}$ with $\mathcal L_\phi h =e^{P(\phi)}h$.
   \end{itemize}
\end{prop}

Therefore, the condition that the transfer operator $\mathcal L_\phi$ has an eigenfunction for  maximal eigenvalue
$\lambda=e^{P(\phi)}$ is equivalent to the existence of an equilibrium state for $\phi$ on $X_+$ having this eigenfunction as Radon-Nikodym derivative with respect to some Gibbs state $\nu\in \G(X_+,\phi)$.

We should note that the approach based on Proposition \ref{existence of eigenfunction-criteria} has already been used in at least two particular cases: for  Dyson potentials by A. Johansson, A. Öberg, M. Pollicott \cite{JOP2023}, and for product-type potentials by L. Cioletti, M. Denker, A. Lopes, M. Stadlbauer \cite{CDLS2017}. Let us now recall these results. The Dyson potential $\phi^{D}:\{-1,1\}^{\Z_+}\to \R$  is given by 
\begin{equation}\label{Dyson potential}
    \phi^{D}(x):=\beta\sum_{n=1}^{\infty}\frac{x_0x_n}{n^{\alpha}},
    \end{equation}
where $\beta\geq 0$ is the inverse temperature and $\alpha>1$ is the model  parameter. 
The Dyson potential $\phi^{D}$ originates from the standard Dyson interaction $\Phi$, by means of (\ref{phi to Phi association}), 
\begin{equation}\label{Dyson interaction with zero m.f.}
    \Phi_{\Lambda}(\omega):=\begin{cases} 
      -\frac{ \beta\omega_i\omega_j}{|i-j|^{\alpha}}, & \text{ if }\Lambda=\{i,j\}\subset\Z, \ i\neq j; \\
      0, & \text{otherwise}.
   \end{cases}
\end{equation}

If $\alpha>2$, 
$\phi^{D}$ has summable variations, and, hence, the Ruelle-Walters theorem applies, and the transfer operator
has a unique positive continuous eigenfunction \cite{Walters1978}. For $\alpha\in (1,2]$, in complete analogy to the classical (whole-line)
Dyson model on $\{-1,1\}^{\Z}$,  phase transitions occur \cite{JOP2019}, and no general
result in Dynamical Systems applies. The main result of \cite{JOP2023} reads
\begin{thm}\cite{JOP2023}\label{JOP2023}
   For $\alpha\in \big(\frac{3}{2},2\big]$ and all sufficiently small $\beta\in[0,+\infty)$ there exists a positive continuous eigenfunction of the Perron-Frobenius transfer operator $\La_{\phi^{D}}$. 
\end{thm}
\begin{remark}
    Theorem \ref{JOP2023} holds for all $\beta<\beta_c^1$, where $\beta_c^1$ is the critical value for a certain long-range Bernoulli percolation.
    For further details see Section \ref{final remarks} and \cite{JOP2023}.
\end{remark}

In \cite{CDLS2017}, the authors introduced a product-type potential $\phi^{P}:\{-1,1\}^{\Z_+}\to \R$, with
\begin{equation}\label{product-type potential}
    \phi^{P}(x):=\beta\sum_{n=1}^{\infty}\frac{x_n}{n^{\alpha}},
    \end{equation}
where again $\beta\geq 0$ and $\alpha>1$. As above, $\phi^{P}$ has a  summable variation for $\alpha>2$, and thus the standard theory applies \cite{Walters1978}. 
For each $\beta$ and $\alpha>1$, there is a unique Gibbs state $\nu$ which 
has the product form $\nu=\prod_{n=0}^\infty \lambda_n$, with $\lambda_n(1)=p_n=\frac{\exp{(\beta\sum_{i=1}^ni^{-\alpha})}}{2\cosh{(\beta\sum_{i=1}^ni^{-\alpha})}}$,
while the unique equilibrium state
  for $\phi^P$ is the Bernoulli measure  $\mu$ on $\{-1,1\}^{\Z_+}$ with $\mu([1]_0)=\lim\limits_{n\to\infty}p_n=\frac{e^{\beta\zeta(\alpha)}}{2\cosh (\beta\zeta(\alpha))}$.

\begin{thm}\cite{CDLS2017}\label{product-type result by CDLS}
    \begin{itemize}        \item[(i)] For $\alpha>3/2$, 
    the equilibrium state 
        $\mu^P$ is absolutely continuous with respect to the Gibbs state $\nu^P$,
 and thus the Perron-Frobenius transfer operator $\La_{\phi^{pt}}$ has an eigenfunction $h^P:=\frac{d\mu^{P}}{d\nu^{P}}\in L^1(X_+,\nu^{P})$. The density $h^P$ is continuous for $\alpha>2$, and essentially discontinuous if $\alpha\le 2$.
        
        \item[(ii)] If $1<\alpha\leq 3/2$, then $\mu^{P}$ and $\nu^{P}$ are singular measures, and therefore,
      the transfer operator $\La_{\phi^{P}}$ does not have an eigenfunction in $L^1(X_+, \nu^{P})$. 
    \end{itemize}
\end{thm}

\section{Basic notions III:  Equivalence of interactions.  Dobrushin uniqueness condition}\label{Basic notions III:  Equivalence of interactions.  Dobrushin uniqueness condition}
 
 \subsection{Equivalent specifications}
We start this section with the notion of equivalent Gibbsian specifications; we follow closely the book by Georgii 
\cite[Chapter 7]{Georgii-book}.

\begin{defin} Two Gibbsian specifications $\gamma, \tilde \gamma$ are called equivalent,
denoted by $\gamma  \simeq \tilde \gamma $ if there exists a constant $C>1$
such that
$$
c^{-1} \gamma_\Lambda(A \mid \cdot) \leq \tilde \gamma_\Lambda(A \mid \cdot) \leq c \gamma_\Lambda(A \mid \cdot)
$$
for all $ \Lambda \in \mathbb{Z}$ and $A \in \mathcal F$.
\end{defin}

In particular, if  $\Phi$ and $\Psi$ are both   UAC interactions, and $\gamma=\gamma^\Phi$, 
$\tilde \gamma=\gamma^\Psi$ are the corresponding Gibbsian specifications, then 
$\gamma^\Phi  \simeq  \gamma^{\Psi} $ if
\begin{equation}\label{cond:hamilequiv}
    \sup _{\Lambda \Subset \mathbb{Z}}\left\|H_{\Lambda}^\Phi-H_{\Lambda}^\Psi\right\|<\infty.
\end{equation}
The sufficient condition (\ref{cond:hamilequiv}) certainly holds if the collection of sets
$$
\left\{\Lambda \Subset \mathbb{Z}: \Phi_\Lambda(\cdot) \neq \Psi_{\Lambda}(\cdot)\right\} 
$$
is finite.  In this case, we say that $\Psi$ is a \textit{\textbf{finite perturbation}} of $\Phi$, and vice versa.

The following theorem summarizes the properties of sets of Gibbs measures of two equivalent specifications.

\begin{thm}\label{equivalent specs}\cite[Theorem 7.3]{Georgii-book}
Let $\gamma$ and $\tilde{\gamma}$ be two equivalent specifications. 
Then $\mathscr{G}(\gamma) \neq \emptyset$ if and only if $\mathscr{G}(\tilde{\gamma}) \neq \emptyset$, and in this case there is an affine bijection $\mu \leftrightarrow \tilde{\mu}$ between $\mathscr{G}(\gamma)$ and $\mathscr{G}(\tilde{\gamma})$ 
such that $\mu=\tilde{\mu}$ on $\mathscr{T}$.
 In particular, $|\operatorname{ex} \mathscr{G}(\gamma)|=|\operatorname{ex} \mathscr{G}(\tilde{\gamma})|$.
\end{thm}

\subsection{Dobrushin uniqueness condition and its corollaries}\label{DUC subsection}
This is one of the most general criteria for the uniqueness of Gibbs states.  We discuss it in the framework of a general countable set $\bL$ of sites and a configuration space $\Omega:=E^{\bL}$. 
Consider a uniformly absolutely convergent (UAC) interaction $\Phi=\{\Phi_\Lambda(\cdot): \Lambda\Subset\bL\}$ on $\Omega$ and let $\gamma^{\Phi}$ be the corresponding Gibbsian specification. For any sites $i,j\in\bL$, define 
$$
C(\gamma^{\Phi})_{i,j}:=\sup_{\eta_{\bL\setminus\{j\}}=\overline{\eta}_{\bL\setminus\{j\}}}||\gamma_{\{i\}}^{\Phi}(\cdot|\eta)-\gamma^{\Phi}_{\{i\}}(\cdot|\overline{\eta})||_{\infty},
$$ 
where $||\cdot||_{\infty}$ is the supremum norm on $\mathcal{M}(\Omega)$ defined by $||\tau||_{\infty}:=\sup_{B\in\mathcal{B}(\Omega)}|\tau(B)|$  for any finite signed Borel measure $\tau$.
The infinite matrix  $C(\gamma^{\Phi}):=(C(\gamma^{\Phi})_{i,j})_{i,j\in\bL}$ is called the \emph{Dobrushin interdependence matrix}.
\medskip

\begin{defin}
The specification $\gamma^{\Phi}$ satisfies the \textbf{\textit{Dobrushin uniqueness (contraction) condition}} if
\begin{equation}\label{DUC for specification}
    c(\gamma^{\Phi}):=\sup_{i\in\bL}\sum_{j\in\bL}C(\gamma^{\Phi})_{i,j}<1. 
\end{equation}
\end{defin}

The Dobrushin uniqueness condition admits a slightly stronger --- and easy to check--- form:
\begin{prop}\cite[Proposition 8.8]{Georgii-book}\label{DUC for interactions}
Let $\bL$ be any countable set, and suppose  $\Phi=\{\Phi_\Lambda(\cdot): \Lambda\Subset\bL\}$ is such that
\begin{equation}\label{DUC for general interactions}
\bar c(\Phi):=
\frac 12 \sup_{i\in\bL}\sum_{\Lambda\ni i} (|\Lambda|-1)\,\delta(\Phi_\Lambda)<1,
\end{equation}
where  $\delta(f):=\sup\{|f(\xi)-f(\eta)|:\xi,\eta\in E^{\bL}\}$ is the variation of $f:E^{\bL}\to\R$. 
Then $\gamma^{\Phi}$ satisfies the Dobrushin uniqueness condition.
\end{prop}

The proof of Proposition \ref{DUC for interactions} boils down to 
showing that for all $i,j$, $i\ne j$,
$$C(\gamma^{\Phi})_{ij}\le 
\frac 12 
\sum_{\Lambda\ni\{i,j\}}
\delta({\Phi_\Lambda})=:\bar C(\Phi)_{ij},$$
and hence
$c(\gamma^{\Phi})=\sup_i \sum_{j} C(\gamma^\Phi)_{ij}\le\sup_i \sum_{j} \bar C(\Phi)_{ij}=\bar c(\Phi)$. 
Notice that the non-negative matrix $\bar C(\Phi):=(\bar C(\Phi)_{i,j})_{i,j\in\bL}$  is symmetric.  

Note that by transitioning to condition (\ref{DUC for general interactions}), we are reinforcing the primary condition of this paper, which is the Dobrushin uniqueness condition (\ref{DUC for specification}), but with a particular purpose.
The condition (\ref{DUC for general interactions}) is stable under a \textit{perturbation} of the underlying model/interaction $\Phi$.
Indeed, let $\Psi=\{\Psi_V\}_{V\Subset\bL}$ be an interaction such that for $V\Subset\bL$, either $\Psi_V=\Phi_V$ or $\Psi_V=0$, then it is straightforward to check that $\bar c(\Psi)\leq \bar c(\Phi)$. 
Thus $\Psi$ inherits the Dobrushin uniqueness condition from $\Phi$ as long as $\Phi$ satisfies (\ref{DUC for general interactions}).

The crucial property of the Dobrushin uniqueness condition is that it provides the uniqueness of the compatible probability measures with the specification $\gamma^\Phi$. 
In fact, we have the following theorem.
\begin{thm}\cite[Chapter 8]{Georgii-book} If  $\gamma^{\Phi}$ satisfies the Dobrushin uniqueness condition (\ref{DUC for specification}), then
$\left|\G(\gamma^{\Phi})\right|\le 1$.
\end{thm}

If $E$ is compact, as is the case in this paper, there is always at least one Gibbs state; hence, the inequality becomes an equality.

The validity of Dobrushin's criterion yields two important properties of the unique Gibbs state: 
concentration inequalities and explicit bounds on the decay of correlations.
\\
The first property involves the coefficient $\bar{c}({\Phi})$, and it provides the tail bounds for the unique Gibbs measure.
Let
\begin{equation}
\delta_kF:=\sup\bigl\{F(\xi)-F(\eta):\xi_j=\eta_j,\; j\in\bL\setminus\{k\}\bigr\}
\end{equation}
denote the variation of a local function $F:\Omega\to\R$ at a site $k\in\bL$, where $\Omega=E^{\bL}$.

\begin{thm}{\cites{K2003}}\label{Kulske-GCB for general potential tails}
    Suppose $\Phi$ is a UAC interaction satisfying (\ref{DUC for general interactions}) and let 
    $\mu_{\Phi}$ be its unique Gibbs measure. 
    Set 
    \begin{equation}
    D:=\frac{4}{(1-\bar{c}(\Phi))^2}\;.
    \end{equation}
    Then, for all $t>0$ and every continuous function $F$ on $\Omega$, one has
    \begin{equation}\label{Kulske tail bounds}
        \mu_{\Phi}\Big(\Big\{\omega\in\Omega: \; F(\omega)-\int_{\Omega} F d\mu_{\Phi}\;\geq\; t\Big\}\Big)\;\leq\; e^{-\frac{2t^2}{D||\underline{\delta}(F)||_2^2}},
    \end{equation}
    where $||\underline{\delta}(F)||^2_{2}:=\sum_{k\in\bL} (\delta_kF)^2$.
\end{thm}

{
It is a well-known fact that (\ref{Kulske tail bounds}) implies that $F$ is \textit{sub-Gaussian} and $\mu_\Phi$ has the \textit{moment concentration bounds} as stated in the following theorem(\cite[Proposition 2.5.2]{V2018}). 
\begin{thm}\cite{V2018}\label{GCB and MCB}
Assume that a probability measure $\mu_\Phi$ satisfies (\ref{Kulske tail bounds}) with a constant $D=D(\mu_\Phi)>0$. 
Then:
    \begin{itemize}
        \item[(i)] $\mu_{\Phi}$ satisfies a \textit{\textbf{Gaussian Concentration Bound}} with  the constant $D$, i.e., for any continuous function $F$ on $\Omega=E^{\bL}$, one has
    \begin{equation}\label{GCB ineq}
     \int_{\Omega}e^{F-\int_{\Omega}Fd\mu_\Phi}d\mu_\Phi\leq e^{D||\underline{\delta}(F)||^2_{2}}.
    \end{equation}
    \item[(ii)] for all $m\in\N$ and any continuous function $F$ on $\Omega$, one has
   \begin{equation}
       \int_{\Omega} \Big|F-\int_{\Omega}Fd\mu_{\Phi}\Big|^m d\mu_{\Phi}
       \leq \Big(\frac{D||\underline \delta (F)||_2^2}{2} \Big)^{\frac{m}{2}} m \Gamma\Big(\frac{m}{2}\Big),
   \end{equation}
   where $\Gamma$ is Euler's gamma function. 
    \end{itemize}
\end{thm}

We present the second property in the particular setup $\bL=\Z$, and it involves the $\Z\times\Z$ matrix
\begin{equation}
D(\gamma^{\Phi})=\sum_{n=0}^{\infty} C(\gamma^{\Phi})^n\;.
\end{equation}
 The sum of the  $\Z\times\Z$ matrices in the right-hand side converges due to the Dobrushin condition (\ref{DUC for specification}).

\begin{prop}\cites{Follmer1982, Georgii-book} Consider a UAC interaction $\Phi$ on $X=E^\Z$.
\begin{itemize}
    \item[(i)] Assume the specification  $\gamma^{\Phi}$ satisfies the Dobrushin condition (\ref{DUC for specification}) and let $\mu$ be its unique Gibbs measure.
    Then, for all $f,g\in C(X)$ and $i\in\Z$, 
    \begin{equation}\label{Follmer's result for each term}
    \Big|\cov_{\mu}(f,g\circ S^i)\Big| 
    \leq\frac{1}{4}\sum_{k,j\in\Z} D(\gamma^{\Phi})_{jk}\cdot  \delta_kf\cdot\delta_{j-i}g.
    \end{equation}
    \item[(ii)] { Suppose $\Phi$ satisfies (\ref{DUC for general interactions}) and define the  non-negative symmetric $\Z\times \Z$-matrix by
    \begin{equation}
    \bar D(\Phi):=\sum_{n\geq 0}\bar C(\Phi)^n\;.
    \end{equation}
    Then,
    \begin{equation}\label{Follmer's result-property of D}
   \sup_{i\in\Z}\sum_{j\in\Z} \bar D(\Phi)_{ij}\leq \frac{1}{1-\bar c(\Phi)}. 
\end{equation}
}
\end{itemize}

\end{prop}

}

\bigskip

\section{Main results I:  From half-line  to whole-line specifications and measures}\label{novel class of potentials}

\subsection{From half-line to whole-line specifications}
The main results of our paper are grouped into two parts.  In the first part, we consider a
 half-line potential $\phi:X_+\to\R$ satisfying minor
technical assumptions and we identify the natural translation-invariant specification on $X=E^{\Z}$ which \emph{extends} $\gl^\phi$.
In the second part we provide sufficient conditions for the 
Gibbs/equilibrium state
on $X_+$ to be absolutely continuous/equivalent with respect to the corresponding half-line Gibbs state $\nu$.
Then Proposition 3.1 allows us to conclude that the transfer operator admits an eigenfunction.\\
First, we consider the issue of whether there is a whole-line specification naturally associated with a half-line specification $\gl^\phi$. It turns out that under a very mild condition on $\phi$ we obtain an affirmative answer.

\begin{defin}
 We say that a continuous potential $\phi:X_+\to\R$ satisfies the \textbf{extensibility} condition if for all $a_0,b_0\in E$ the sequence 
 $$
 F_n^{a_0,b_0}(x):=S_{n+1}\phi(x_{-n}^{-1}b_0x_1^{\infty})-S_{n+1}\phi(x_{-n}^{-1}a_0x_1^{\infty})
 =\sum_{i=0}^n \bigl(\phi(x_{-i}^{-1}b_0x_1^{\infty})-\phi(x_{-i}^{-1}a_0x_1^{\infty})\bigr)
 $$
 converges uniformly on $x\in X$ as $n\to\infty$.
\end{defin}

This condition has first appeared in \cite{BFV2019} in connection
to a related question of whether $g$-measures are also Gibbs. 
In terms of the extensibility condition, we can reformulate the main result of \cite{BFV2019} as follows.
\begin{thm}\cite{BFV2019}
   A $g-$measure $\mu$ (the natural extension) for a normalized function $g\in \G(X_+)$ becomes a Gibbs measure on $X$ if and only if $\log g$ satisfies the extensibility condition.  
\end{thm}

In \cite{BFV2019}, the authors identified several sufficient conditions for  $\phi$ to satisfy the extensibility condition, such as the Walters condition and the so-called Good Future condition. 
Note that a function $\phi\in C(X_+,\R)$ satisfies the Walters condition and the Good Future condition if $\lim_{p\to\infty} \sup_{n\geq 1}v_{n+p}(S_{n}\phi)=0$ and $\sum_{k=1}^\infty \delta_k(\phi)<\infty$, respectively.  Here  
\begin{equation}\label{eq:variations}
    v_k(\varphi):=\sup_{x_0^{k-1}=y_0^{k-1}}|\varphi(x)-\varphi(y)|
\end{equation}
is the $k^{\text{th}}-$variation (or the oscillation in the volume $[0,k-1]$) of a function $\varphi:X_+\to\R$
and
\begin{equation}
\label{dobcoef}
\delta_k(\varphi) \equiv \sup_{x \in X_+,a_k,b_k 
\in E} \left| \varphi(x_0^{k-1} a_k x_{k+1}^\infty ) - \varphi(x_0^{k-1} b_k x_{k+1}^\infty)\right|.
\end{equation}
is the oscillation of $\varphi$ at the site $k$.
Interesting examples of potentials satisfying the extensibility condition are the Dyson potential (\ref{Dyson potential}) and the product-type potential (\ref{product-type potential}).
Note that both potentials are in the Walters class if and only if $\alpha>2$. However, both satisfy the Good Future condition, and thus the extensibility condition as well, for all admissible values of the parameter $\alpha$, since $\delta_n(\phi^{D})=\mathcal{O}(n^{-\alpha})$ and $\delta_n(\phi^{P})=\mathcal{O}(n^{-\alpha})$. 

Gibbsianness for potentials satisfying the extensibility condition stems from the fact that they lead to a natural whole-line translation-invariant specification. Indeed, if $\phi\in C(X_+,\R)$ satisfies the extensibility condition the limits 
\begin{equation}\label{whole-line specification for extensibility function}
   \overset{\rightleftharpoons}{\gamma}^{\phi}_{\{i\}}(\sigma_i|\omega_{\{i\}^c}):=\lim_{p\to\infty} \frac{e^{S_{i+p+1}\phi(\omega_{-p}^{i-1}\sigma_i\omega_{i+1}^{\infty})}}{\sum_{\Bar{\omega}_i}e^{S_{i+p+1}\phi(\omega_{-p}^{i-1}\Bar{\omega_i}\omega_{i+1}^{\infty})}}
\end{equation}
are well defined for all  $i\in\Z$.  In turn, they lead to a full specification.

\begin{prop}\label{existence of two sided specifications for extensibility potentials}
There is a unique translation-invariant quasilocal non-null specification $\overset{\rightleftharpoons}{\gamma}^{\phi}$ on $X$ such that the functions  (\ref{whole-line specification for extensibility function}) become its single-site densities. 
\end{prop}
\begin{proof}
By Proposition \ref{important theorem about single-site densities}, it is sufficient to check that $\{\overset{\rightleftharpoons}{\gamma}^{\phi}_{\{i\}}\}_{i\in \Z}$ satisfy conditions (i) and (ii) of Definition \ref{def single site densit}. Clearly, $\{\overset{\rightleftharpoons}{\gamma}^{\phi}_{\{i\}}\}_{i\in \Z}$ satisfies the first condition
In order to check the second condition, consider arbitrary $i,j\in\Z$ with $i<j$, and $p\in \N$ such that $i,j\gg-p$. Then 
 a straightforward computation {shows} that
\begin{multline*}
    \frac{e^{S_{i+p+1}\phi(\omega_{-p}^{i-1}\alpha_i\omega_{i+1}^{j-1}\alpha_j\omega_{j+1}^{\infty})}}{\sum_{\beta_{i,j}}e^{S_{j+p+1}\phi(\omega_{-p}^{i-1}\beta_i\omega_{i+1}^{j-1}\beta_j\omega_{j+1}^{\infty})+S_{i+p+1}\phi(\omega_{-p}^{i-1}\beta_i\omega_{i+1}^{j-1}\alpha_j\omega_{j+1}^{\infty})-S_{j+p+1}\phi(\omega_{-p}^{i-1}\beta_i\omega_{i+1}^{j-1}\alpha_j\omega_{j+1}^{\infty}) }}=\\
    =\frac{e^{S_{j+p+1}\phi(\omega_{-p}^{i-1}\alpha_i\omega_{i+1}^{j-1}\alpha_j\omega_{j+1}^{\infty})}}{\sum_{\beta_{i,j}}e^{S_{i+p+1}\phi(\omega_{-p}^{i-1}\beta_i\omega_{i+1}^{j-1}\beta_j\omega_{j+1}^{\infty})+S_{j+p+1}\phi(\omega_{-p}^{i-1}\alpha_i\omega_{i+1}^{j-1}\beta_j\omega_{j+1}^{\infty})-S_{i+p+1}\phi(\omega_{-p}^{i-1}\alpha_i\omega_{i+1}^{j-1}\beta_j\omega_{j+1}^{\infty}) }}.
\end{multline*}
The limit $p\to\infty$ of this identity yields  condition (ii). 

Finally, note that for all $i\in \Z$,
$$
\overset{\rightleftharpoons}{\gamma}^{\phi}_{\{0\}}((S^i\sigma)_0|S^i(\omega)_{\{0\}^c})=\lim_{p\to\infty}\frac{\exp{(S_{p+1}\phi(\omega_{i-p}^{i-1}\sigma_i\omega_{i+1}^{\infty}))}}{\sum_{\Bar{\sigma}_i}\exp{(S_{p+1}\phi(\omega_{i-p}^{i-1}\Bar{\sigma}_i\omega_{i+1}^{\infty}))}}=\overset{\rightleftharpoons}{\gamma}^{\phi}_{\{i\}}(\sigma_i|\omega_{\{i\}^c})
$$
and hence  $\overset{\rightleftharpoons}{\gamma}^{\phi}$ is a translation-invariant specification.
\end{proof}

\begin{example}
    Let us illustrate the construction in the previous proof with the Dyson potential $\phi$ defined in (\ref{Dyson potential}). 
    In this case a straightforward computation shows that for all $p\in\N$, $\beta\geq 0$, $\sigma_0\in\{-1,1\}$ and $\omega\in \{-1,1\}^{\Z}$,
    $$
    \frac {\exp\left( {S_{p+1}\phi( \omega_{-p}^{-1}\sigma_0\omega_1^\infty)}\right)}
    { \exp\left( {S_{p+1}\phi( \omega_{-p}^{-1}\bar \sigma_0\omega_1^\infty)}\right)}
    =
    \frac {\exp\left( \beta\sum_{k=-p}^{-1} \frac {\sigma_0\omega_k}{|k|^\alpha}+\beta\sum_{k=1}^{+\infty} \frac {\sigma_0\omega_k}{|k|^\alpha}\right)}
    {\exp\left(\beta\sum_{k=-p}^{-1} \frac {\bar \sigma_0\omega_k}{|k|^\alpha}+\beta\sum_{k=1}^{+\infty} \frac {\bar \sigma_0\omega_k}{|k|^\alpha}\right)}.
    $$
    Thus 
    $$
    \frac {\exp\left( {S_{p+1}\phi( \omega_{-p}^{-1}\sigma_0\omega_1^\infty)}\right)}
    {  {\exp\left( {S_{p+1}\phi( \omega_{-p}^{-1}\sigma_0\omega_1^\infty)}\right)}
    +\exp\left( {S_{p+1}\phi( \omega_{-p}^{-1}\bar \sigma_0\omega_1^\infty)}\right)}\to 
    \gamma_{\{0\}}(\sigma_0|\omega_{-\infty}^{-1},\omega_1^\infty)
    $$
    as $p\to\infty$, where 
    $$
    \gamma_{\{0\}}(\sigma_0|\omega_{-\infty}^{-1},\omega_1^\infty)
    =
    \frac {\exp\left(\beta\sum_{k=1}^{+\infty} \frac {\sigma_0(\omega_k+\omega_{-k}}{k^\alpha}\right)}
    {\exp\left(\beta\sum_{k=1}^{+\infty} \frac {\sigma_0(\omega_k+\omega_{-k}}{k^\alpha}\right)
    +
    \exp\left(\beta\sum_{k=1}^{+\infty} \frac {\bar\sigma_0(\omega_k+\omega_{-k}}{k^\alpha}\right)}.
    $$
\end{example}

\subsection{Whole-line Gibbs measures as limits of half-line ones}\label{from half-line Gibbs to whole-line Gibbs}
Once the correspondence between $\gl^{\phi}$ and $\overset{\rightleftharpoons}{\gamma}^{\phi}$ is established, it is interesting to understand the relation between $\G(X_+,\gl^{\phi})\Big(=\G(X_+,\phi)\Big)$ and $\G(X,\overset{\rightleftharpoons}{\gamma}^{\phi})$.

Suppose $\nu$ is a Gibbs probability  measure on $X_+$ for some potential $\phi$ satisfying the extensibility condition, 
and consider an arbitrary measure  $\rho$ on $X_{-}=E^{-\mathbb N}$(e.g.  a uniform Bernoulli measure),
and consider the measure $\mu_0=\rho\times\nu$ on $X=X_{-}\times X_+$.

For any $n\ge 0$, let
$$
\mu_n = \mu_0\circ S^{-n},
$$
This is a sequence of Borel probability measures $\{\mu_n\}_{n\ge 0}$ on the compact 
metric space $X$, and hence it has weak$^*$-converging subsequences $\{\mu_{n_k}\}$.  It turns out that their limits are bona fide whole-line Gibbs measures.

\begin{manualtheorem}{A}\label{convergence of shifted measures} 
Suppose $\phi$ satisfies the extensibility condition and
$\nu\in\G(X_+,\phi)$. Consider $\mu_0:=\rho\times\nu$, where $\rho$ is uniform Bernoulli measure on $X_-$.
Assume that a subsequence $\{\mu_{n_k} =\mu_0\circ S^{-n_k}\}_{k\ge 0}$ converges to a probability measure $\mu$ in the weak$^*$ topology as $k\to\infty$.
Then for ${
\mu}$-almost all $x$:

$$
\mu(x_0|x_{-\infty}^{-1},x_{1}^\infty)
=\overset{\rightleftharpoons}{\gamma}^{\phi}_{\{0\}}(x_0|x_{-\infty}^{-1},x_{1}^\infty)=\lim_{n\to\infty}
\frac {\exp\left( S_{n+1}\phi(x_{-n}^{-1} x_0 x_{1}^\infty)\right)}
 {\sum_{\bar a_0}\exp\left( S_{n+1}\phi(x_{-n}^{-1} \bar a_0 x_{1}^\infty)\right)}.
$$

Hence $\mu$ is a whole-line Gibbs measure for the whole-line specification defined by the kernels
$\overset{\rightleftharpoons}{\gamma}^{\phi}_{\{0\}}$.
\end{manualtheorem}

Now we turn to translation invariant measures.

\begin{manualtheorem}{B}\label{convergence of averaged measures} 
Assume the conditions of Theorem \ref{convergence of shifted measures}, let $\mu$ be a weak$^*$ limit point of the following sequence of measures 
$$\Tilde{\mu}_n:=\frac{1}{n}\sum_{i=0}^{n-1}\mu_0\circ S^{-i}.$$ 
Then $\mu$ is a translation-invariant Gibbs measure for $\overset{\rightleftharpoons}{\gamma}^{\phi}$, i.e., $\mu\in \G_{S}(\overset{\rightleftharpoons}{\gamma}^{\phi}):=\G(\overset{\rightleftharpoons}{\gamma}^{\phi})\cap \mathcal{M}_1(X,S)$.
\end{manualtheorem}
We shall give the proofs of Theorem \ref{convergence of shifted measures} and Theorem \ref{convergence of averaged measures} in Section \ref{proofs of the main results}.

\section{Main results II: Eigenfunctions of transfer operators and absolute continuity of Gibbs measures}\label{Eigenfunctions of transfer operators and absolute continuity of Gibbs measures}

As mentioned in the introduction,
we can show the existence of an eigenfunction
if we can show that the  equilibrium state $\mu$ is absolutely continuous with respect to
the half-line non-translation invariant 
Gibbs measure $\nu$. 
This approach has recently been used by Johansson, \"Oberg and Pollicott \cite{JOP2023} to show Theorem \ref{JOP2023}. We also use this approach to prove the following theorem.

\begin{manualtheorem}{C}\label{main theorem about the Dyson model} 
Let $\phi$ be the Dyson potential (\ref{Dyson potential}). Suppose $\alpha>1$ and $\beta>0$ is sufficiently small.
Then, 
\begin{itemize}
    \item[(i)] the half-line  Dyson model $\phi$ on $X_+=\{\pm 1\}^{\Z_+}$ admits a unique equilibrium state $\mu_+$ and Gibbs state $\nu$;
    \item[(ii)] for all $\alpha>1$, $\mu_+$ is equivalent to $\nu$, i.e., $\mu_+\ll\nu$ and $\nu\ll\mu_+$, and thus the Perron-Frobenius transfer operator $\La_{\phi}$ has an eigenfunction in $L^1(X_+,\nu)$;  
    \item[(iii)] furthermore, if $\alpha>\frac{3}{2}$, 
    there exists a continuous version
    of the Radon-Nikodym density $\frac{d\mu_+}{d\nu}$, and thus the Perron-Frobenius transfer operator $\La_{\phi}$ has a continuous eigenfunction. 
\end{itemize}
\end{manualtheorem}
\begin{remark}
    In fact, we prove Theorem \ref{main theorem about the Dyson model} under the Dobrushin uniqueness condition (c.f. (\ref{DUC for specification})-(\ref{DUC for general interactions})). 
    Therefore, the parameter $\beta\geq0$ should be sufficiently small so that the corresponding whole-line Dyson interaction on $X=\{\pm 1\}^{\Z}$ satisfies the Dobrushin uniqueness condition (\ref{DUC for general interactions}).
\end{remark}

We will show a much more general statement (Theorem \ref{sufficient cond for UI} below), from which one can easily obtain the first two parts of Theorem \ref{main theorem about the Dyson model}.
Our method of proof of statement (iii) differs from that of \cite{JOP2023}, and can be generalized to other models. 
The proof of Theorem \ref{main theorem about the Dyson model} will be presented in Section \ref{proofs of the main results}.

\subsection{From $C(X_+)$ to the space of interactions on $X$}\label{From $C(X_+)$ to the space of interactions on $X$}

It is clear from Section \ref{novel class of potentials} that every (half-line) translation-invariant UAC interaction $\gPhi$ (i.e., $\gPhi_{\Lambda}\circ S=\gPhi_{\Lambda+1}$ 
) on $X_+$, yields a potential $\phi\in C(X_+)$ with $\phi=-\sum_{0\in V\Subset \Z_+}\gPhi_V$, such that $\gamma^{\phi}=\gamma^{\gPhi}$.  

In the opposite direction, it is known that every potential $\phi$ yields an equivalent interaction $\gPhi$ on $X_+$ \cite{Rue} which in its turn can be extended by translations to an interaction $\Phi$ on $X$. This extension, however, is only guaranteed to belong to the so-called $\mathscr{B}_{0}(X)$ class ----i.e., is such that $\sum_{0\in V\Subset\Z_+} \frac{1}{|V|}||\Phi_V||_{\infty}<\infty$--- but it may fail to be UAC. 

The determination of necessary and sufficient conditions for $\phi$ to yield a UAC potential $\Phi$ on $X$ is an open problem that we do not address here.  Rather, in the sequel we determine a
class of potentials $\phi$ that admit a translation invariant UAC interaction $\Phi$ on $\Z$ satisfying 
\begin{equation}\label{standing assumption about phi}
    \phi=-\sum_{0\in V\Subset\Z_+}\Phi_V.
\end{equation}
The following proposition shows that potentials in this class satisfy the extensibility condition.
\begin{prop}\label{phi-to-Phi_the form of corresponding whole-line spec}
    Let $\phi\in C(X_+)$ such that there exists a translation invariant UAC interaction $\Phi$ on $\Z$ satisfying 
    $\phi=-\sum_{0\in V\Subset\Z_+}\Phi_V$. 
    Then $\phi$ satisfies the $extensibility$ condition. 
    Furthermore, $\glr^{\phi}=\gamma^{\Phi}$.
\end{prop}
\begin{proof}
Note that 
$$
S_{n+1}\phi=-\sum\limits_{V\cap [0,n]\neq\emptyset;\; V\Subset\Z_+}\Phi_V.
$$ 
Thus since the interaction is translation-invariant ($\Phi\circ S^{-n}=\Phi_{V-n}$)
$$
S_{n+1}\phi\circ S^{-n}=-\sum\limits_{V\cap [-n,0]\neq \emptyset;\; V\Subset[-n,\infty]}\Phi_V.
$$
(The LHS of the above equation should be considered on $X$, otherwise, $S^{-n}$ is not defined on $X_+$.)

Pick any $\xi,\eta\in X$, such that $\xi_j=\eta_j$ if $j\neq 0$. 
Then
\begin{equation*}
    S_{n+1}\phi(\xi_{-n}^{\infty})-S_{n+1}\phi(\eta_{-n}^{\infty})
    =\sum\limits_{V\cap [-n,0]\neq \emptyset;\; V\Subset[-n,\infty]}(\Phi_V(\eta)-\Phi_V(\xi))=\sum\limits_{0\in V\Subset[-n,\infty]}(\Phi_V(\eta)-\Phi_V(\xi)) .
\end{equation*}

Since $\Phi$ is UAC,  $\sum\limits_{0\in V\Subset[-n,\infty]}[\Phi_V(\xi)-\Phi_V(\eta)]$ converges uniformly to $H^{\Phi}_{\{0\}}(\xi)-H^{\Phi}_{\{0\}}(\eta)=\sum_{0\in V\Subset\Z}(\Phi_V(\xi)-\Phi_V(\eta))$ as $n\to\infty$. 
Thus it is also clear that $\glr^{\phi}=\gamma^{\Phi}$.
\end{proof}

\begin{example}\label{Association Dyson potential and interaction}
    Consider the Dyson potential $\phi$ (c.f. (\ref{Dyson potential})) with the state space $E=\{-1,1\}$, and the standard (whole-line) Dyson interaction $\Phi$ on $\Z$ (c.f. (\ref{Dyson interaction with zero m.f.})). 
Then it is easy to see that 
$
\phi=-\sum_{0\in\Lambda \Subset \Z_+}\Phi_{\Lambda}.
$
\end{example}

\subsection{Decoupling across the origin}
To every
UAC interaction $\Phi$ on $\Z$ we can associate a sequence  ${\Psi}^{(k)}$ of interactions obtained by removing all bonds linking $-\N$ and $\Z_+$ and adding them one at a time. Formally, we consider the family
$$
\mathcal A=\bigl\{\Lambda\Subset Z:\, \min(\Lambda)<0,\ \max \Lambda\ge 0\bigr\}.
$$
indexed according to some arbitrary order: $\mathcal A=\{\Lambda_1,\Lambda_2,\ldots\}.$
Then define, for each $k\in\Z_+$, 
\begin{equation}\label{construction of intermediate interactions}
    {\Psi}^{(k)}_\Lambda = \begin{cases}  0,&\quad \Lambda\in\{\Lambda_i: i\geq k+1\},\\
\Phi_{\Lambda}, &\quad \text{ otherwise}.
\end{cases}\;
\end{equation}
In particular, $\Psi^{(0)}$ has no interaction between the left and right half lines.  
Clearly, all the constructed interactions are UAC and in addition have the following properties: 
\begin{remark}\label{convergence of intermediate interactions and specs }
    \begin{itemize}
        \item[1)] Every $\Psi^{(k)}$ in (\ref{construction of intermediate interactions}) is a local (finite) perturbation of $\Psi^{(0)}$. Moreover, the sequence $\Psi^{(k)}$ tends to $\Phi$
        as $k\to\infty$, in the sense that 
        $\Psi^{(k)}_{\Lambda}\rightrightarrows\Phi_{\Lambda}$ 
        for all $\Lambda\Subset \Z$.
        \item[2)] For any finite volume $V$, \
    
        $$||H^{\Psi^{(k)}}_V-H^{\Phi}_V||_{\infty}\leq \sum\limits_{\substack{\Lambda_j\cap V\neq \emptyset \\j\geq k}}||\Phi_{\Lambda_j}||_{\infty}\xrightarrow[k\to\infty]{}0. 
        $$
        \item[3)] The specifications $\gamma^{\Psi^{(k)}}$ converge to $\gamma^{\Phi}$ as $k\to\infty$.  More precisely, for all $B\in\F$ and $V\Subset\Z$,
        $$
\gamma^{\Psi^{(k)}}_V(B|\omega)\xrightarrow[k\to\infty] 
    {}\gamma^{\Phi}_V(B|\omega)\; \text{ 
        uniformly in the boundary conditions } 
        \omega\in X.
        $$
        \item[4)] In addition, if $\nu^{(k)}$ is a Gibbs measure for $\Psi^{(k)}$, then by Lemma 4.2, any weak$^*$-limit point, $\mu$ of the sequence $\{\nu^{(k)}\}_{k\geq 0}$ is a Gibbs measure for the potential $\Phi$.

    \end{itemize}
\end{remark}

Another important observation applies to the interaction $\Psi^{(0)}$.  As it is
  constructed from $\Phi$ by removing all interaction between $-\N$ and $\Z_+$,  the corresponding specification $\gamma^{\Psi^{(0)}}$ is of product type \cite[Example 7.18]{Georgii-book}:
 $\gamma^{\Psi^{(0)}}=\gamma^{\Phi^-}\times\gamma^{\Phi^+}$, where $\Phi^-$ and $\Phi^+$ are the restrictions of $\Phi$ respectively to the negative and positive half-lines. Thus, 
 the extreme Gibbs measures also factorize \cite[Example 7.18]{Georgii-book}:

\begin{equation}\label{extreme Gibbs measures for Psi^0}
    \text{ex}\;\G(\gamma^{\Psi^{(0)}})
    =\{\nu_-\times\nu_+: \nu_-\in \text{ex}\;\G(X_-,\gamma^{\Phi^-}),\;\; \nu_+\in \text{ex}\;\G(X_+,\gamma^{\Phi^+}) \}\;.
\end{equation}

\subsection{Absolute continuity}
If $\Psi^{(0)}$ does not exhibit phase transitions, then by Theorem \ref{equivalent specs}, all the interactions $\Psi^{(k)}$, $k\geq 1$, do not exhibit phase transitions as well.
Let us denote by $\nu^{(k)}$ the unique Gibbs state for $\Psi^{(k)}$.

\begin{thm}\label{equiv of intermediate Gibbs meas}
   If the interaction $\Phi$ satisfies the Dobrushin uniqueness  criterion (\ref{DUC for general interactions}), so do all interactions $\Psi^{(k)}$
    and, furthermore, $\nu^{(k)}$ and $\nu^{(0)}$ are equivalent with
    \begin{equation}\label{densities with respect to the product of half-line Gibbs measures}
   \frac{d\nu^{(k)}}{d\nu^{(0)}}=\frac{e^{-\sum_{i=1}^k\Phi_{\Lambda_i}}}{\int_{X}e^{-\sum_{i=1}^k\Phi_{\Lambda_i}}d\nu^{(0)}} 
   \quad,\quad   \frac{d\nu^{(0)}}{d\nu^{(k)}}=\frac{e^{\sum_{i=1}^k\Phi_{\Lambda_i}}}{\int_{X}e^{\sum_{i=1}^k\Phi_{\Lambda_i}}d\nu^{(k)}}\;.
\end{equation}
   
\end{thm}
\begin{proof}
The hereditary character of Dobrushin's criterion follows from the obvious fact  
 that, for all $k\in\Z_+$, $\bar c(\Psi^{(k)})\leq \bar c(\Phi)$. 

The proof of (\ref{densities with respect to the product of half-line Gibbs measures}) follows, telescopically, from the partial Radon-Nikodym derivatives
\begin{equation}\label{densities-the recursion}
    \frac{d\nu^{(k)}}{d\nu^{(k-1)}} = \frac {e^{-\Phi_{\Lambda_k}}}
 { \int e^{-\Phi_{\Lambda_{k}}} d\nu^{(k-1)}}
 \quad, \quad
  \frac{d\nu^{(k-1)}}{d\nu^{(k)}} = \frac {e^{\Phi_{\Lambda_k}}}
 { \int e^{\Phi_{\Lambda_{k}}} d\nu^{(k)}
 }\;,
\end{equation}
which is a particular case of the following elementary lemma.
\end{proof}

\begin{lem}\label{important elementary lemma}
    Let ${\Psi}$ be a UAC interaction that does not exhibit phase transitions.  
    Consider a perturbed interaction of the form $\bar{\Psi}=\Psi+\FP$ with $\FP$ a finite interaction supported on $A\Subset \Z$.  
    If $\{{\nu}\}=\G({\Psi})$ and $\Bar{\nu}\in\G(\Bar{\Psi})$, then $\Bar{\nu}\ll {\nu}$, and 
    \begin{equation}\label{eq:rfr.1}
        \frac{d\bar{\nu}}{d{\nu}}= \frac{e^{-H_A^\FP}}{\int_{X}e^{-H_A^\FP}d{\nu}}.
    \end{equation}
\end{lem}
\begin{proof}
    First, note that by Theorem \ref{equivalent specs}, the model $\bar \Psi$ does not exhibit phase transitions since it is a finite perturbation of $\Psi$. Thus $\{\Bar{\nu}\}=\G(\Bar{\Psi})$.  By uniqueness of the Gibbs state, we have that, for both interactions, the Gibbs measures are achieved through limits 
    \begin{equation}\label{eq:rfr.0}
{\gamma}_{\Lambda}(\cdot)\stackrel{\ast}{\xrightharpoonup[\Lambda\uparrow\Z]{}}{\nu}\; \mbox{and}\ \bar{\gamma}_{\Lambda}(\cdot)\stackrel{\ast}{\xrightharpoonup[\Lambda\uparrow\Z]{}}\bar{\nu}
\end{equation}
where ${\gamma}:=\gamma^{{\Psi}}$ and 
    $\bar{\gamma}:=\gamma^{\bar{\Psi}}$ are the Gibbsian kernels corresponding to the interactions ${\Psi}$ and $\bar{\Psi}$ and free boundary conditions (that is, considering only bonds within $\Lambda$). Consider any $V\subset\Lambda\Subset\Z$.  As
 $H_\Lambda^{\bar\Psi}=H_A^\FP + H_\Lambda^\Psi$, for any cylindrical event $[\sigma_V]$, we have that 
 \begin{equation*}
\bar{\gamma}_\Lambda(\sigma_V)
            \;=\;\frac{\sum_{\xi_{\Lambda}} \mathbf{1}_{\sigma_V}\, e^{-H^\FP_A(\xi_A)}\,e^{-H_{\Lambda}^{{\Psi}}(\xi_{\Lambda})}}{\sum_{\xi_{\Lambda}}e^{-H^\FP_A(\xi_A)}\,e^{-H_{\Lambda}^{{\Psi}}(\xi_{\Lambda})}}.
\end{equation*}
Dividing top and bottom by $\sum_{\xi_\Lambda}e^{-H_{\Lambda}^{{\Psi}}(\xi_{\Lambda})}$ we get
\begin{equation*}
\bar{\gamma}_\Lambda(\sigma_V)
            \;=\;\frac{\gamma_\Lambda\bigl (\mathbf{1}_{\sigma_V}\, e^{-H^\FP_A}\bigr)}{\gamma_\Lambda\bigl (e^{-H^\FP_A}\bigr)}.
\end{equation*}
    Taking limits over $\Lambda$ and using (\ref{eq:rfr.0}) we obtain
    \begin{equation*}
        \bar{\nu}(\mathbf{1}_{\sigma_V})\;=\; \nu(\mathbf{1}_{\sigma_V}\,f)
    \end{equation*}
where $f$ is, precisely, the right-hand side of (\ref{eq:rfr.1}).  This concludes the proof because cylindrical events uniquely determine the measures.  
    \end{proof}
Now we turn to the restrictions $\nu^{(k)}_+$ of measures $\nu^{(k)}$ to the half-line $\Z_+$, i.e., for $B\in \F_+$, $\nu^{(k)}_+(B):=\nu^{(k)}(X_-\times B)$. 
Note that for the cylindrical sets $[\sigma_{\Lambda}]$, $\Lambda\Subset\Z_+$, one has $\nu^{(k)}_+([\sigma_{\Lambda}])=\nu^{(k)}([\sigma_{\Lambda}])$. 
Similarly, one can define the restrictions to the left half-line $-\N$. 
If the model $\Psi^{(0)}$ does not exhibit phase transitions, then $\nu^{(0)}$ is a product measure, i.e., $\nu^{(0)}=\nu^{(0)}_-\times \nu^{(0)}_+$ (c.f.  (\ref{extreme Gibbs measures for Psi^0})), where $\nu^{(0)}_-$ is the unique Gibbs measure for the interaction $\Phi^-$ on $X_-$, and $\nu^{(0)}_+$ is the unique Gibbs measure for $\Phi^+$ on $X_+$. If this is the case, then one can compute the Radon-Nikodym density $f^{(k)}_+:=\frac{d\nu^{(k)}_+}{d\nu^{(0)}_+}$, in fact, for all $\sigma\in X_+$, $k\in\Z_+$,
\begin{equation}\label{Radon-Nikodym densities on Z_+}
    f^{(k)}_{+}(\sigma) 
    =\frac {\int_{X_-} e^{-\sum_{j=1}^{k} \Phi_{\Lambda_j} (\xi,\sigma)} \nu_-^{(0)}(d\xi)}
    {\int_{X_+}\int_{X_-} e^{-\sum_{j=1}^{k} \Phi_{\Lambda_j}(\xi,\zeta)}\nu^{(0)}_-(d\xi)\nu^{(0)}_+(d\zeta)}.
\end{equation}

By Theorem \ref{equiv of intermediate Gibbs meas}, all the measures $\nu^{(k)},\; k\geq 0$ are equivalent to $\nu^{(0)}$.
However, it is not clear whether the weak$^*$-limit points of the sequence $\{\nu^{(k)}\}_{k\in\Z_+}$ are absolutely continuous with respect to $\nu^{(0)}$ or not.
The following theorem provides sufficient conditions.  
\begin{manualtheorem}{D}\label{UI and Dunford-Pettis}
    Assume that $\Phi$ satisfies the Dobrushin uniqueness condition (\ref{DUC for general interactions}). 
    Suppose the family $\{f^{(k)}\}_{k\in\N}$ is uniformly integrable in $L^1(\nu^{(0)})$. 
    Then the weak$^*$ limit point of the sequence $\{\nu^{(k)}\}$ is a Gibbs measure for $\Phi$ and absolutely continuous with respect to $\nu^{(0)}$.
\end{manualtheorem}
The next theorem is the main theorem of this section and it provides sufficient conditions for uniform integrability of the family $\{f^{(k)}\}_{k\in\N}$, and thus absolute continuity of the weak$^*$ limit points with respect to $\nu^{(0)}$.
\begin{manualtheorem}{E}\label{sufficient cond for UI}
    Assume the following
    \begin{itemize}
        \item[1)] the interaction $\Phi$ satisfies the Dobrushin uniqueness condition (\ref{DUC for general interactions});
        \item[2)] $\sum_{k=1}^{\infty}||\underline{\delta}(\Phi_{\Lambda_k})||_2^2<\infty$;
        \item[3)] $\sum_{k=1}^{\infty}\rho_k<\infty$, where $\rho_k:=\sup\limits_{n\in\N}\left|\int_{X}\Phi_{\Lambda_k}d{\nu}^{(n)}\right|$.
    \end{itemize}
    Then $\Psi^{(0)}$ does not exhibit phase transitions, and $\{f^{(k)}:k\in\N\}$ is uniformly integrable in $L^1(\nu^{(0)})$.
\end{manualtheorem}
\begin{remark}
Note that the third, summability condition in Theorem \ref{sufficient cond for UI} 
is important.
To illustrate this, consider the product-type potential $\phi$ defined by (\ref{product-type potential}). 
One can readily verify that the potential $\phi$ coincides with the half-line mean energy at $0$ --- that is, $\phi=-\sum_{0\in V\Subset \Z_+}\Phi_V$--- for the
  translation-invariant UAC interaction:
$$
\Phi_{\Lambda} (\omega)= \begin{cases}  -\frac{\beta\omega_j}{|i-j|^{\alpha}},&\quad \Lambda=\{i,j\}\subset \Z, \;j>i,\\
0, &\quad \text{ otherwise }.
\end{cases}\;
$$
This interaction $\Phi$, which is not spin-flip invariant, satisfies the first and second conditions of Theorem \ref{sufficient cond for UI} for $\beta$ sufficiently small. 
In fact, for all $\beta\geq 0$, both Dobrushin interdependence matrices $C(\gamma^{\Phi})$ and $C(\gamma^{\phi})$ are zero matrices, thus, both specifications (not the interactions) satisfy the condition (\ref{DUC for specification}) for all $\beta$'s. 
However, $\Phi$ does not satisfy the last condition of the theorem. 
Indeed, as the unique Gibbs measure $\mu$ for $\Phi$ is Bernoulli with 
\[
\mu([1]_0)=\frac{e^{\beta\zeta(\alpha)}}{2\cosh (\beta\zeta(\alpha))},
\]
one has 
\[
\int_X\sigma_0d\mu=\tanh(\beta\zeta(\alpha))>0.
\]
Therefore,  
$$
\rho_{-i,j}\geq \Big|\int_X \Phi_{\{-i,j\}}d\mu\Big|
=\frac{\beta}{(i+j)^{\alpha}}\tanh(\beta\zeta(\alpha)),
$$
for all $i\in\N,j\in\Z_+$, and thus the sum $\sum_{i=1}^{\infty}\sum_{j=0}^{\infty}\rho_{-i,j}$ diverges.
\end{remark}

The combination of Theorems \ref{sufficient cond for UI} and \ref{UI and 
Dunford-Pettis} implies that if $\Phi$ satisfies the conditions of the former, 
the unique Gibbs measure $\mu\in\G(\Phi)$ is absolutely continuous with 
respect to $\nu^{(0)}$. 
Note that, as a consequence, Theorem  \ref{sufficient cond for UI} can not be 
true for the previous interaction if $1<\alpha\leq 3/2$.  Indeed, the second 
part of Theorem \ref{product-type result by CDLS} directly implies that  the 
measures $\mu$ and $\nu^{(0)}$ are singular for those values of $\alpha$.

It is natural to ask whether, reciprocally,  $\nu^{(0)}$ is also absolutely continuous with respect to $\mu$. 
The answer is affirmative modulo conditions comparable to those of Theorem \ref{sufficient cond for UI}.  The argument resorts also to a sequence of interactions but these time obtained by removing one by one the bonds in the volumes in $\A=\{\Lambda\Subset\Z: \min \Lambda<0, \max \Lambda\geq 0\}$ instead of adding them to $\Psi^{(0)}$.  More precisely, at each step $k\in\Z_+$, we construct a new interaction $\Phi^{(k)}$ as follows:
$$
{\Phi}^{(k)}_\Lambda = \begin{cases}  0,&\quad \Lambda\in\{\Lambda_i: 1\leq i\leq k\},\\
\Phi_{\Lambda}, &\quad \text{ otherwise}.
\end{cases}\;
$$
Then as previously, for all $k\in\Z_+$, the matrix $\bar C(\Phi^{(k)})$ is dominated by $\bar C(\Phi)$, thus if $\Phi$ satisfies the Dobrushin uniqueness condition (\ref{DUC for general interactions}) so do all  the interactions $\Phi^{(k)}$. 
Furthermore, Remark \ref{convergence of intermediate interactions and specs } still remains valid interchanging  $\Phi$ and $\Psi^{(0)}$.
Thus the sequence $\mu^{(k)}$ of unique Gibbs measures for each $\Phi^{(k)}$ converges in weak$^*$ sense to a Gibbs measure $\mu$ for $\Psi^{(0)}$. 
In addition, by Lemma \ref{important elementary lemma}, all the measures $\mu^{(k)}$ are equivalent to $\mu$ and their Radon-Nikodym derivatives are given by 
\[
f^{(k)}_{\mu}:=\frac{d\mu^{(k)}}{d\mu}=\frac{e^{\sum_{i=1}^{k}\Phi_{\Lambda_i}}}{\int_X e^{\sum_{i=1}^{k}\Phi_{\Lambda_i}}d\mu}.
\]
To conclude, we present  the following analogue of Theorem \ref{sufficient cond for UI}
which can be proven in a similar way. 
\begin{manualtheorem}{F}\label{equivalence of mu and nu^(0)}
    Assume that
    \begin{itemize}
        \item[1)] the interaction $\Phi$ satisfies (\ref{DUC for general interactions});
        \item[2)] $\sum_{k=1}^{\infty}||\underline{\delta}(\Phi_{\Lambda_k})||_2^2<\infty$;
        \item[3)] $\sum_{k=1}^{\infty}\rho_k^{\mu}<\infty$, where $\rho_k^{\mu}:=\sup\limits_{n\in\N}\Big|\int_{X}\Phi_{\Lambda_k}d{\mu}^{(n)}\Big|$.
    \end{itemize}
    Then $\Psi^{(0)}$ does not exhibit phase transitions, and $\{f^{(k)}_{\mu}:k\in\N\}$ is uniformly integrable in $L^1(\mu)$. In particular, $\nu^{(0)}$ is absolutely continuous with respect to $\mu$.
\end{manualtheorem}

We postpone the proofs of Theorem \ref{UI and Dunford-Pettis} and Theorem \ref{sufficient cond for UI} until Section \ref{proofs of the main results}.

\section{Application: Dyson model}\label{Applications: Dyson model}
This was our motivating example.  Let us 
recall that the Dyson potential $\phi$ is defined on the half-line  configuration space $X_+=\{-1,+1\}^{\Z_+}$ as 
$$
\phi( x) = \beta\sum_{n=1}^\infty \frac {x_0x_n}{n^\alpha},\quad x=(x_0,x_1,\ldots)\in X_+,
$$
for some $\alpha>1$.
As mentioned in Example \ref{Association Dyson potential and interaction}, this potential is 
related to the whole-line Dyson model $\Phi$, defined in (\ref{Dyson interaction with zero m.f.}), by
$\phi=-\sum_{0\in V\Subset\Z_+}\Phi_V$,
and thus
$\overset{\rightleftharpoons}{\gamma}^{\phi}=\gamma^{\Phi}$.
Hence 
by Theorem \ref{convergence of shifted measures}, for all $\nu\in\G(\phi
)$, 
any weak* accumulation point of the sequence $\{\nu\circ S^{-n}\}$ 
is a restriction of some  $\mu\in \G(\Phi)$ to $\Z_+$.

The Dyson potential satisfies the Good Future condition because its oscillations (\ref{dobcoef}) $\delta_n(\phi) = \mathcal O(n^{-\alpha})$ are summable.  As a consequence, it satisfies the extensibility condition. 
For $\alpha>2$,
furthermore, its variations $v_{n}(S_n\phi):=\sup\{S_n\phi(x)-S_n\phi(y):x_{0}^{n-1}=y_0^{n-1}\}$ are also summable and,  therefore, the standard theory applies \cite{Walters1978}.  
The case of $\alpha\in (1,2]$ is significantly more subtle, and its theory is less developed. 
Its most recent advance is (Theorem \ref{JOP2023} above), obtained by { Johansson, \"Oberg} and Pollicott \cite{JOP2023} 
through 
the random-cluster representation for the whole-line Dyson model $\Phi$ (c.f. (\ref{Dyson interaction with zero m.f.})).

Now, we recall some important properties of the Dyson model and its phase diagram.
It is clear that the interaction $\Phi$ is translation and spin-flip invariant. 
One of the most interesting properties of the associated Gibbs measures is that all the measures in $\G(\Phi)$ are translation-invariant for all the values of the parameters $\alpha$ and $\beta$.
However, this statement is not true as regards spin-flip invariance. 
In fact, there is only one spin-flip invariant Gibbs measure in $\G(\Phi)$ which is the cause of the phase transitions for high $\beta$'s. 
By applying FKG inequalities, it can be shown that the (weak$^*$) limits $\mu^+:=\lim_{\Lambda \uparrow \mathbb{Z}}\gamma^{\Phi}_{\Lambda}(\cdot|+)$ and $\mu^-:=\lim_{\Lambda \uparrow \mathbb{Z}}\gamma^{\Phi}_{\Lambda}(\cdot|-)$ exist, and both are extremal. 
In fact, $\mu^-$ and $\mu^{+}$ are the only extremal elements of $\G(\Phi)$ since they stochastically dominate all other Gibbs measures of the model i.e. $\int_X f d\mu^-\leq \int_X f d\mu\leq \int_X f d\mu^+$ for any Gibbs state $\mu\in \G(\Phi)$ and for any non-decreasing function $f\in C(X)$. The phase diagram of this model is in many ways similar to the phase diagram of the two-dimensional nearest-neighbour Ising model, in fact, we have the following result.

\begin{thm}\cites{Dyson1969,Ruelle1972,FS1982}\label{phase diagram of whole-line Dyson model}
    For all $\alpha\in(1,2]$, there exists a critical temperature $\beta_c(\alpha)\in (0,+\infty)$ such that there is no 
    phase transition for all $\beta\in [0,\beta_c(\alpha))$ (i.e. $\mu^{+}=\mu^{-}$), and there is a phase transition for all $\beta\in(\beta_c(\alpha),+\infty)$ (i.e. $\mu^{+}\neq\mu^{-}$). 
    Furthermore, for all the values of $\beta$, $\G(\Phi)=[\mu^{-};\mu^{+}]$.
\end{thm}
Note that if $\alpha>2$, then for all the values of the inverse tempereature $\beta>0$, $\Phi$  does not exhibit phase transitions, i.e., $|\G(\Phi)|=1$.

In \cite{JOP2019}, the authors proved that the phase diagram of the half-line Dyson model $\phi$ is similar to the phase diagram of the whole-line Dyson model. In fact, they showed that for all $\alpha\in(1,2]$, there exists 
$\beta_c^+$, such that the for all $\beta\in (0,\beta^+_c)$, there exists a unique half-line Gibbs state for $\phi$, and for $\beta>\beta^+_c$, there exist multiple Gibbs states.
The authors also conjectured that the 
critical values $\beta_c^+$ and $\beta_c$ of
the half and whole-line Dyson models are, in fact,  equal $\beta^+_{c}=\beta_c$.

\section{Proofs of the main results}\label{proofs of the main results}
\subsection{Proofs of Theorems \ref{convergence of shifted measures} and \ref{convergence of averaged measures}}
We shall use the following two simple lemmas in the proof of Theorem \ref{convergence of shifted measures}.

\begin{lem}\cite[Remark 5.10]{Georgii-book}\label{Gibbsianity of shifted measure}
Let $\gamma$ be a specification on $X$ 
and 
 $\mu\in \G(\gamma)$. Then $\mu\circ S^{-1}$ is consistent with the specification $\gamma^{(1)}$, where 
$$\gamma^{(1)}_{\Lambda}(B|\omega):=\gamma_{\Lambda+1}(S^{-1}(B)|S^{-1}(\omega)), \; \forall B\in\F, \; \omega\in X, \; \forall\Lambda\Subset \Z .$$
\end{lem}
\begin{lem}\cite[Theorem 4.17]{Georgii-book}\label{convergence of spec implies convergence of Gibbs meas}
Suppose $\gamma$ and $\gamma^{(n)}, \; n\geq 1$ are specifications on $X=E^{\Z
}$.
Assume that $\gamma^{(n)}$ converge uniformly to $\gamma$ as $n\to\infty$, in the sense that  for all $\Lambda\Subset\Z$ and all $\sigma\in X$, $$\gamma_{\Lambda}(\sigma_{\Lambda}|\omega_{\Lambda^c})=\lim_{n\to\infty}{\gamma}^{(n)}_{\Lambda}(\sigma_{\Lambda}|\omega_{\Lambda^c})$$ uniformly in the boundary condition $\omega\in X$. 
Take $\mu^{(n)}\in \G(\gamma^{(n)})$, and assume that the sequence $\mu^{(n)}$ converges to some $\mu\in\mathcal{M}_1(X)$ in the weak$^*$ topology. 
Then $\mu\in\G(\gamma)$.
\end{lem}

\begin{proof}[Proof of Theorem \ref{convergence of shifted measures}.]
Let us consider a family of functions $\gamma_{\{n\}}^{(0)}:E\times E^{\Z\setminus\{n\}}\to (0,1)$ given by:
\begin{equation}\label{auxiliary specification}
    \gamma_{\{n\}}^{(0)}(a_n|x_{-\infty}^{n-1},x_{n+1}^{\infty}):=\begin{cases} 
     \gl^\phi_{\{n\}}(a_n|x_{-\infty}^{n-1},x_{n+1}^{\infty}), & n\geq 0; \\
      1/|E|, & n<0
   \end{cases}
\end{equation}
with
\begin{equation}
\gl^\phi_{\{n\}}(a_n|x_{-\infty}^{n-1},x_{n+1}^{\infty})\;=\; \frac{e^{S_{n+1}\phi(x_0^{n-1}a_n x_{n+1}^{\infty})}}{\sum_{\Bar{a}_n}e^{S_{n+1}\phi(x_0^{n-1}\Bar{a}_n x_{n+1}^{\infty})}}
\end{equation}
are the single-site kernels for the half-line  specification $\gl^\phi$.
[To simplify we adopt the convention $x_m^n=\emptyset$ if $n<m$.] It can be easily checked that it is a family of single-site densities of a specification, therefore, by Proposition \ref{important theorem about single-site densities} there is a non-null quasilocal specification $\gamma^{(0)}$ on $(X,\mathcal{F})$ having $\{\gamma_{\{n\}}^{(0)}\}_{n\in\Z}$ as its single-site densities. 
Note that $\gamma^{(0)}$ is not 
 translation-invariant.
We claim that  $\mu_0\in\G({\gamma}^{(0)})$.  To show this, it is enough, by Proposition \ref{important theorem about single-site densities}, to check that $\mu_0\bigl( \gamma_{\{n\}}^{(0)}(F)\bigr)=\mu_0(F)$ for each local cylindrical $F$.  That is, we must show that
\begin{equation}\label{eq:consistency}
[\rho\times\nu]\bigl( \gamma_{\{n\}}^{(0)}(\mathds{1}_{b_k^r})\bigr)\;=\; [\rho\times\nu]\bigl( \mathds{1}_{b_k^r}\bigr)
\end{equation}
for all integer $k\le r$ and $n$ and all $b_k^r\in E^{r-k+1}$.  A quick inspection shows that the only non-trivial case is $k\le n \le r$, $n\ge 0$.  
In this case
\begin{equation}
[\rho\times\nu]\bigl( \gamma_{\{n\}}^{(0)}(\mathds{1}_{b_k^r})\bigr)\;=\; \rho\bigl( \mathds{1}_{b_{k\wedge(-1)}^0}\bigr)\,
\nu\bigl( \gl_{\{n\}}(\mathds{1}_{b_{k\vee 0}^r})\bigr)\;,
\end{equation} 
and, as  $\nu$ is a Gibbs measure for $\phi$,  we conclude that
\begin{equation}
[\rho\times\nu]\bigl( \gamma_{\{n\}}^{(0)}(\mathds{1}_{b_k^r})\bigr)\;=\; \rho\bigl( \mathds{1}_{b_{k\wedge(-1)}^0}\bigr)\,
\nu\bigl( \mathds{1}_{b_{k\vee 0}^r}\bigr)\;,
\end{equation} 
proving (\ref{eq:consistency}).

The proof is concluded by invoking the previous lemmas.  By Lemma \ref{Gibbsianity of shifted measure}, 
for all $p\in\N$,
the measure $\mu_p=\mu_0\circ S^{-p}$ is consistent with $\gamma^{(p)}$, where 
$$\gamma^{(p)}_{\Lambda}(B|\omega):=\gamma^{(0)}_{\Lambda+p}(S^{-p}(B)|S^{-p}(\omega)), \; \forall B\in\F, \; \omega\in X, \; \forall\Lambda\Subset \Z.$$ 
Note that the single-site density functions of $\gamma^{(p)}$ can be calculated explicitly, namely,
for all $\sigma, \omega\in X$,
\begin{equation}\label{p-auxiliary specification}
    \gamma_{\{i\}}^{(p)}(\sigma_i|\omega_{\{i\}^c}):=\begin{cases} 
      \frac{e^{S_{i+p+1}\phi(\omega_{-p}^{i-1}\sigma_i \omega_{i+1}^{\infty})}}{\sum_{\Bar{\omega}_i}e^{S_{i+p+1}\phi(\omega_{-p}^{i-1}\Bar{\omega}_i \omega_{i+1}^{\infty})}}, & i\geq -p; \\
      {1}/{|E|}, & i<-p.
   \end{cases}
\end{equation}
Thus $\overset{\rightleftharpoons}{\gamma}^{\phi}$ is the uniform limit of the sequence of specifications $\{\gamma^{(p)}\}_p$.
Thus by Lemma \ref{convergence of spec implies convergence of Gibbs meas}, we obtain that $\mu\in\G(\overset{\rightleftharpoons}{\gamma}^{\phi})$, and hence we conclude the statement of the theorem. 
\end{proof}

The following lemma will be useful in the proof of Theorem \ref{convergence of averaged measures}.
\begin{lem}\label{Corollary of Theorem 6.1.}
Under the conditions of Theorem \ref{convergence of shifted measures},  for all cylindrical sets $C\subset X$, and all volumes $\Lambda\Subset \Z$, 
\[
\lim_{n\to\infty}\int_{X}\Big[\overset{\rightleftharpoons}{\gamma}^{\phi}_{\Lambda}(C|x)-\mathds{1}_{C}(x)\Big]\mu_n(dx)=0\;.
\]
\end{lem}
\begin{proof}
We will show that 
\begin{equation}\label{eq:liminf}
\liminf_{n\to\infty}\int_{X}\Big[\overset{\rightleftharpoons}{\gamma}^{\phi}_{\Lambda}(C|x)-\mathds{1}_{C}(x)\Big]\mu_n(dx)=0\;;
\end{equation}
the analogous result for the $\limsup$ can be shown similarly.
Take any subsequence $\{n_k\}$ such that
$$\lim_{k\to\infty}\int_{X}\Big[\overset{\rightleftharpoons}{\gamma}^{\phi}_{\Lambda}(C|x)-\mathds{1}_{C}(x)\Big]\mu_{n_k}(dx)=\liminf_{n\to\infty}\int_{X}\Big[\overset{\rightleftharpoons}{\gamma}^{\phi}_{\Lambda}(C|x)-\mathds{1}_{C}(x)\Big]\mu_n(dx).$$ 
By compactness of $\mathcal{M}_1(X)$ the subsequence $\{\mu_{n_k}\}_k$ converges in the weak$^*$-topology.  If  $\mu$ is its limit, then $\mu\in\G(\overset{\rightleftharpoons}{\gamma}^{\phi})$ by Theorem \ref{convergence of shifted measures}.  Thus, 
\[
\lim_{k\to\infty}\int_{X}\Big[\overset{\rightleftharpoons}{\gamma}^{\phi}_{\Lambda}(C|x)-\mathds{1}_{C}(x)\Big]\mu_{n_k}(dx)=
\int_{X}\Big[\overset{\rightleftharpoons}{\gamma}^{\phi}_{\Lambda}(C|x)-\mathds{1}_{C}(x)\Big]\mu(dx)
=0\;. 
\]
\end{proof}

\begin{proof}[Proof of Theorem \ref{convergence of averaged measures}.]
The proof of translation-invariance of $\mu$ is standard. Thus it is enough to check the consistency, i.e., $\mu\in \G(\overset{\rightleftharpoons}{\gamma}^{\phi})$.
Let $\mu=\lim_{k}\Tilde{\mu}_{n_k}$, and take any cylindrical event $C$. Then, the weak convergence implies that
\begin{equation}\label{Eq1 in Theorem 6.2.}
\int_X\overset{\rightleftharpoons}{\gamma}^{\phi}_{\Lambda}(C|x)\Tilde{\mu}_{n_k}(dx)\xrightarrow[k\to\infty]{}(\mu\overset{\rightleftharpoons}{\gamma}^{\phi}_{\Lambda})(C).    
\end{equation}
On the other hand,  the Stolz-Cesaro theorem and Lemma \ref{Corollary of Theorem 6.1.} yield that
\[
\frac{1}{n_k}\sum_{i=0}^{n_k-1}\int_{X}\Big[\overset{\rightleftharpoons}{\gamma}^{\phi}_{\Lambda}(C|x)-\mathds{1}_{C}(x) \Big]\mu_i(dx)\xrightarrow[k\to\infty]{}0\;.
\]
Thus,
$$\int_X\overset{\rightleftharpoons}{\gamma}^{\phi}_{\Lambda}(C|x)\Tilde{\mu}_{n_k}(dx)=\frac{1}{n_k}\sum_{i=0}^{n_k-1}\int_{X}\Big[\overset{\rightleftharpoons}{\gamma}^{\phi}_{\Lambda}(C|x)-\mathds{1}_{C}(x) \Big]\mu_i(dx)+\Tilde{\mu}_{n_k}(C)\xrightarrow[k\to\infty]{}\mu(C),$$
\end{proof}

\subsection{Proofs of Theorem \ref{UI and Dunford-Pettis}, \ref{sufficient cond for UI} and \ref{main theorem about the Dyson model}}
\begin{proof}[Proof of Theorem \ref{UI and Dunford-Pettis}.]
    By Theorem \ref{equiv of intermediate Gibbs meas}, each measure $\nu^{(k)}$ is absolutely continuous with respect to $\nu^{(0)}$ with $f^{(k)}=\frac{d\nu^{(k)}}{d\nu^{(0)}}$. 
Furthermore, by Theorem \ref{equivalent specs}, $\nu^{(k)}$ is the unique Gibbs measure for $\Psi^{(k)}$ for all $k\geq 0$. 
Let $\mu^*$ be a weak$^*$ limit of a subsequence $\{\nu^{(k_s)}\}_{s\in\N}$. 
By Lemma \ref{convergence of spec implies convergence of Gibbs meas}, $\mu^{*}\in\G(\Phi)$ (c.f. Remark \ref{convergence of intermediate interactions and specs }).   
By the weak star convergence, we have that for all $g_0\in C(X)$,
\begin{equation}\label{weak star convergence of intermediate Gibbs meas}
    \int_{X} g_0d\nu^{(k_s)}=\int_{X}g_0f^{(k_s)}d\nu^{(0)}\xrightarrow[s\to\infty]{}\int_{X}g_0d\mu^*.
\end{equation}
Since the family $\{f^{(k)}:k\in\Z_+\}$ is uniformly integrable, it is relatively weakly compact in $L^1(\nu^{(0)})$ by the Dunford-Pettis theorem. Therefore, there exists a weak limit point $f\in L^1(\nu^{(0)})$ of the sequence $\{f^{(k_s)}\}_{s\in\N}$.  Without loss of generality, assume that $f^{(k_s)}\xrightharpoonup[s\to\infty]{}f$. Thus for all $g\in L^{\infty}(X,\nu^{(0)})$,
\begin{equation}\label{weak convergence of the intermediate densities}
    \int_{X}gf^{(k_s)}d\nu^{(0)}\xrightarrow[s\to\infty]{}\int_{X}gfd\nu^{(0)}.
\end{equation}
By combining (\ref{weak star convergence of intermediate Gibbs meas}) and (\ref{weak convergence of the intermediate densities}), we conclude that for all $g_0\in C(X)$, 
\[
\int_{X}g_0d\mu^*=\int_{X}g_0fd\nu^{(0)}.
\]
\end{proof}

\begin{proof}[Proof of Theorem \ref{sufficient cond for UI}.]
For all $k\in\N$, denote $W_k:=\sum_{i=1}^k\Phi_{\Lambda_i}$.  Our argument relies on two claims:
\medskip

\noindent
{\bf Claim 1:}
\[
\sup_{k\geq 0} \Big|\int_{X}-W_kd\nu^{(k)}\Big| < \infty\;.
\]
\noindent
{\bf Claim 2:}
\begin{equation}\label{UI of numerator}
    \sup_{k\ge 0}\int_{X}e^{-W_k}d\nu^{(0)}<\infty.
\end{equation}

Then, as   
\begin{equation}\label{the sum}
    \int_{X}f^{(k)}\log f^{(k)}d\nu^{(0)}
    =\int_{X}-W_kd\nu^{(k)}-\log \int_{X}e^{-W_k}d\nu^{(0)}\;,
\end{equation} 
these claims imply that 
\begin{equation}
    \sup_{k\geq 0} \int_{X}f^{(k)}\log f^{(k)}d\nu^{(0)}<\infty.\;
\end{equation}
Hence by applying de la Vall\'ee Poussin's theorem to the family $\{f^{(k)}:k\in\N\}$ and to
the function $t\in(0,+\infty)\mapsto t\log t$, one concludes that the family  $\{f^{(k)}:k\in\N\}$ is uniformly integrable in $L^1(\nu^{(0)})$. 
\medskip

The proof of Claim 1 is immediate:
\begin{equation}\label{de la walle poussin-estim for the first part}
    \Big|\int_{X}-W_kd\nu^{(k)}\Big|\leq \sum_{i=1}^k\Big| \int_{X}\Phi_{\Lambda_i}d\nu^{(k)}\Big|\leq \sum_{i=1}^k\rho_i\leq \sum_{i=1}^{\infty}\rho_i <\infty\;.  
\end{equation}

The proof of Claim 2 relies on the Gaussian concentration bounds. Note that 
 for all $k\in\N$,  $\bar c (\Psi^{(0)})\leq \bar c(\Psi^{(k)})\leq \bar c(\Phi)$.
Therefore, the Dobrushin uniqueness condition $\bar c(\Phi)<1$ is inherited by all the intermediate interactions.  Applying the first part of Theorem \ref{GCB and MCB} we see that the (only) measure 
$ \mu\in\G(\Phi)$ and all the intermediate measures $\nu^{(k)}$, $k\geq 0$, satisfy the Gaussian Concentration Bound with the same constant 
$D:=\frac{4}{(1-\bar c(\Phi))^2}$. 
This implies that, for all $k\in\N$, 
\begin{equation}
    \int_{X}e^{-\Phi_{\Lambda_k}}d\nu^{(k-1)}\leq e^{ D||\underline\delta(\Phi_{\Lambda_k})||^2_2} e^{-\int_{X}\Phi_{\Lambda_k}d\nu^{(k-1)}}.
\end{equation}
We combine this inequality with (\ref{densities-the recursion}) to iterate
\begin{eqnarray}
\lefteqn{ \int_{X}e^{-(\Phi_{\Lambda_k}+\Phi_{\Lambda_{k-1}})}d\nu^{(k-2))}}\nonumber\\
&=&
  \int_{X}e^{-\Phi_{\Lambda_{k}}}d\nu^{(k-1))}\, \int_{X}e^{-\Phi_{\Lambda_{k-1}}}d\nu^{(k-2))} \nonumber\\[8pt]
 &\le& 
 e^{ D(||\underline\delta(\Phi_{\Lambda_k})||^2_2+||\underline\delta(\Phi_{\Lambda_{k-1}})||^2_2)} \cdot e^{-(\int_{X}\Phi_{\Lambda_k}d\nu^{(k-1)}+\int_{X}\Phi_{\Lambda_{k-1}}d\nu^{(k-2)})}.
\end{eqnarray}
By induction this yields 
\begin{equation}
    \int_{X}e^{-\sum_{i=1}^k\Phi_{\Lambda_i}}d\nu^{(0)}\leq 
    e^{ D\sum_{i=1}^k||\underline\delta(\Phi_{\Lambda_i})||^2_2}\cdot e^{-\sum_{i=1}^k\int_{X}\Phi_{\Lambda_i}d\nu^{(i-1)}}.
\end{equation}
Thus,  using \eqref{de la walle poussin-estim for the first part}, we have that 
\begin{equation}
    \int_{X}e^{-W_k}d\nu^{(0)}
    \leq
    e^{ D\sum_{i=1}^k||\underline\delta(\Phi_{\Lambda_i})||^2_2}\cdot e^{\sum_{i=1}^k\rho_i}
\end{equation}
for all $k\in\N$, and
\begin{equation}\label{de la walle poussin-estim for the second part}
   \sup_{k\in\N} \int_{X}e^{- W_k}d\nu^{(0)}\leq e^{ D\sum_{i=1}^{\infty}||\underline\delta(\Phi_{\Lambda_i})||^2_2}\cdot e^{\sum_{i=1}^{\infty}\rho_i}<\infty.
\end{equation}
This proves Claim 2 and, hence, concludes the proof of the theorem.
\end{proof}

\begin{proof}[Proof of Theorem \ref{main theorem about the Dyson model}.]\ 
\smallskip

\textit{Part (i):} Its proof is rather straightforward.
We choose $\beta>0$, so that the resulting $\Phi$ satisfies the Dobrushin uniqueness condition $\bar c(\Phi)<1$.  This condition is inherited by $\Psi^{(0)}$, hence both potentials have a unique Gibbs state.  
Furthermore, as direct products of Gibbs measures for the restricted interactions $\Phi^-$ and $\Phi^+$ are Gibbs measures for $\Psi^{(0)}$, neither $\Phi^-$ nor $\Phi^+$ may exhibit phase transitions (c.f. Equation (\ref{extreme Gibbs measures for Psi^0})). 
In particular, $\gamma^{\phi}=\gamma^{\Phi^+}$ admits only one Gibbs state. 
\smallskip

\textit{Part (ii):} We just have to verify the hypotheses of Theorem \ref{sufficient cond for UI}; the absolute continuity follows then from Theorem \ref{UI and Dunford-Pettis}.  We chose $0<\beta\le \beta_{DU}$ with
\[
\beta_{DU}=\big(2\sum_{i=1}^{\infty}\frac{1}{i^{\alpha}}\big)^{-1}=\frac{1}{2}\zeta(\alpha)^{-1}\;.
\]
Hence, by Proposition \ref{DUC for interactions}, for all $\beta\in(0,\beta_{DU})$, $\Phi$ satisfies the Dobrushin uniqueness condition $\bar c(\Phi)<1$.  
\noindent
\textit{Hypothesis 1):} Consequence of the Dobrushin uniqueness criterion.
\smallskip

\noindent
\textit{Hypothesis 2):} For all $i\in\N,j\in\Z_+$ and $k\in\Z$, 
\begin{equation*}
    \delta_k(\Phi_{\{-i,j\}})=\begin{cases} 
      0, & k\notin \{-i,j\}; \\
      \frac{2\beta}{(i+j)^{\alpha}}, & k\in\{-i,j\}.
     \end{cases}
    \end{equation*}
Thus for all $i\in\N,\; j\in\Z_+$, 
Therefore, 
\[
\sum_{\substack{i\in\N,\\j\in\Z_+}}||\underline\delta(\Phi_{\{-i,j\}})||^2_2=\sum_{\substack{i\in\N,\\j\in\Z_+}}\frac{8\beta^2}{(i+j)^{2\alpha}}<\infty\;.
\]

\noindent
\textit{Hypothesis 3):}
We shall use inequalities (\ref{Follmer's result for each term}) and (\ref{Follmer's result-property of D}). We use the notation introduced in Section \ref{DUC subsection}.
Note that for all $k\geq 0$ we have the componentwise domination
\begin{equation}\label{dominance of D matrices}
C(\gamma^{\Psi^{(k)}})_{i,j}
\leq \bar C(\Phi)_{i,j}\quad,\quad D(\gamma^{\Psi^{(k)}})_{i,j}
\leq \bar D(\Phi)_{i,j}\;.
\end{equation}

Thus since $\bar c(\Phi)<1$
all the specifications $\gamma^{\Psi^{(k)}}$ satisfy (\ref{DUC for specification}) (c.f. Section \ref{DUC subsection}).
Applying (\ref{Follmer's result for each term}) to the measures $\nu^{(k)},\; k\geq 0$ and using (\ref{dominance of D matrices}), we see that
\begin{equation}
    \sup_{k\geq 0}\Big | \int_X \sigma_{-m}\cdot\sigma_0 \circ S^n d\nu^{(k)}\Big|
    \leq \frac{1}{4} \sum_{r,j\in\Z} \bar D(\Phi)_{rj} \cdot\delta_r \sigma_{-m} \cdot\delta_{n-j}\sigma_0. 
\end{equation}
Hence for all $m\in\N$ and $n\in\Z_+$,
\begin{equation}
    (m+n)^{\alpha}\rho_{-m,n}
    \leq \frac{1}{4} \sum_{r,j\in\Z}  \bar D(\Phi)_{rj} \cdot\delta_r \sigma_{-m}\cdot\delta_{n-j}\sigma_0.
\end{equation}
Summing and applying inequality (\ref{Follmer's result-property of D}), we obtain
\begin{align}
    \begin{split}
    \sum_{n=0}^{\infty}(m+n)^{\alpha}\rho_{-m,n}
    &\leq \frac{1}{4} \sum_{r,j,n\in\Z} \bar D(\Phi)_{rj} \cdot\delta_r \sigma_{-m}\cdot\delta_{n-j}\sigma_0
    \\&= \frac{1}{2}\sum_{r\in\Z}\Big[\delta_r\sigma_{-m}
    \cdot\sum_{j\in\Z}\bar D(\Phi)_{rj}\Big]
    \\ &\overset{(\ref{Follmer's result-property of D})}{\leq} \frac{1}{1-\bar c (\Phi)}.
    \end{split}
\end{align}
As a consequence, $m^{\alpha}\cdot\sum_{n=0}^{\infty}\rho_{-m,n}\leq (1-\bar c(\Phi))^{-1}$ for all $m\in\N$, which  implies
\begin{equation}
    \sum_{m=1}^{\infty}\sum_{n=0}^{\infty}\rho_{-m,n}
    \leq \frac{1}{1-\bar c (\Phi)}\; \sum_{m=1}^{\infty}\frac{1}{m^{\alpha}}<\infty.
\end{equation}

\smallskip

This concludes the verification of the hypotheses of Theorem \ref{sufficient cond for UI}, which, together with
Theorem \ref{UI and Dunford-Pettis} imply that the unique limit point $\mu$ of the sequence $\{\nu^{(k)}\}_{k\in\Z_+}$ is absolutely continuous with respect to $\nu^{(0)}$. 
As $\nu^{(0)}=(\nu^{(0)})^-\times(\nu^{(0)})^+$ and $(\nu^{(0)})^+$ coincides with the unique Gibbs state $\nu$ of $\phi$, it follows that $\mu\ll \nu^{(0)}$ implies that $\mu^+\ll \nu$.

The proof of $\nu\ll \mu^{+}$ is analogous but using Theorem \ref{equivalence of mu and nu^(0)} instead of Theorem \ref{sufficient cond for UI}.
\medskip

\textit{Part (iii):} It is an application of the Arzela-Ascoli theorem. Let us denote $\xi$ the configurations in $X_-=\{-1,1\}^{-\N}$ and $\sigma$ those in $X_+=\{-1,1\}^{\Z_+}$. For all $N\in\N$, and $(\xi,\sigma)\in X_-\times X_+$, define 
\[
W_N(\xi,\sigma):=\sum_{i=1}^{N}\sum_{j=0}^{N}-\frac{\beta \xi_{-i}\sigma_{j}}{(i+j)^{\alpha}}
\]
and
\begin{equation}\label{the subsequence of the densities on Z_+}
    f^{[N]}_{+}(\sigma):=\frac{d\nu^{[N]}_+}{d\nu^{(0)}_+}(\sigma) =
    \frac {\int_{X_-} e^{-W_N(\xi,\sigma)} \nu_-^{(0)}(d\xi)}
    {\int_{X_+}\int_{X_-} e^{-W_N(\xi,\zeta)}\nu^{(0)}_-(d\xi)\nu^{(0)}_+(d\zeta)}.
\end{equation}
The sequence $\{f^{[N]}_+\}_{N\in\N}$ is a subsequence of the sequence $\{f^{(k)}_+\}_{k\in\Z_+}$ defined in (\ref{Radon-Nikodym densities on Z_+}). 
All $f^{[N]}_+$ are local functions on $X_+$, thus continuous. 
We claim that it is enough to prove that the family $\{f^{[N]}_+:N\in \N\}$ is relatively compact in $C(X_+)$. 
Indeed, if this is true, there exists a function $f_+\in C(X_+)$ and a subsequence $\{f^{[N_k]}_+\}_{k\in\N}$ such that $f^{[N_k]}_+\rightrightarrows f_+$ as $k\to\infty$. 
Thus, by the argument presented in the proof of Theorem \ref{UI and Dunford-Pettis}, $f_+$ is the Radon-Nikodym density of $\mu_+$ with respect to $\nu$,and
\[
\La_{\phi}f_+(x)=e^{P(S,\phi)} f_+(x)
\] 
for all for $x\in X_+$ [in principle, the identity holds for $\nu-$almost all $x\in X_+$, but $\nu$ is fully supported].  
Hence $f_+$ is the continuous eigenfunction of the transfer operator $\La_{\phi}$ corresponding to the largest eigenvalue. 
Note that we can also conclude from this argument that the entire sequence $\{f^{[N]}_+\}_{N\in\N}$ converges in the uniform topology.
\medskip

To conclude, we turn to the proof of the relative compactness of the family $\{f^{[N]}_+:N\in \N\}$.
By the Arzela-Ascoli theorem, it is enough to show that this family is uniformly bounded and equicontinuous. These properties are proven separately.
\smallskip

\noindent
\textit{Uniform boundedness:}
As
\[
\bar c(\Phi^-)=\frac{1}{2}\sup\limits_{i\in -\N}\sum\limits_{i\in V\Subset -\N}(|V|-1)\delta(\Phi_V)\leq \frac{1}{2}\sup\limits_{i\in \Z}\sum\limits_{i\in V\Subset \Z}(|V|-1)\delta(\Phi_V)=\bar c(\Phi)<1\;,
\] 
 the interaction $\Phi^-$ on $-\N$ satisfies the Dobrushin uniqueness condition, therefore, the unique Gibbs measure $\nu_-\in \G(\Phi^-)$ satisfies the Gaussian Concentration Bound (Theorem \ref{GCB and MCB}) with the constant 
 $D=4(1-\bar{c}(\Phi))^{-2}$.
 Fix any $\sigma\in X_+$ and consider  $W_N(\xi,\sigma)$ as a function of $\xi\in X_-$. Clearly, it is a local function, thus by the first part of Theorem \ref{GCB and MCB}, for all $\kappa\in\R$,
\begin{equation}\label{GCB to nu_-}
    \int_{X_-} e^{\kappa W_N(\xi,\sigma)}\nu_-(d\xi)
    \leq e^{D \kappa^2 ||\underline \delta (W_N(\cdot,\sigma))||_2^2} \cdot e^{\kappa \int_{X_-} W_N(\xi,\sigma)\nu_-(d\xi)}.
\end{equation}
First, note that the interaction $\Phi^-$ is invariant under the global spin-flip transformation, therefore, so is the unique Gibbs measure $\nu_-$.
Thus for all $N\in\N$ and $\sigma\in X_+$, $\int_{X_-} W_N(\xi,\sigma)\nu_-(d\xi)=0$.
Second, for all $k\in\N$, 
$$
\delta_{-k}(W_N(\cdot,\sigma))=2\beta \Big| \sum_{j=0}^N \frac{\sigma_j}{(k+j)^{\alpha}}\Big|\leq 2\beta \sum_{j=0}^N \frac{1}{(k+j)^{\alpha}}.
$$
Hence if $\alpha>\frac{3}{2}$, for all $N\in\N$, $\sigma\in X_+$,
\begin{equation}\label{eq:r-variance}
    ||\underline \delta (W_N(\cdot,\sigma))||_2^2\leq 4\beta^2 \sum_{k=1}^{\infty}\Big(\sum_{j=k}^{\infty}\frac{1}{j^{\alpha}} \Big)^2=:4\beta^2 C_1(\alpha)<\infty.
\end{equation}
Then (\ref{GCB to nu_-}) implies that for all $N\in\N$, $\sigma\in X_+$, one has
\begin{equation}\label{UB for num of f^N_+}
    \int_{X_-} e^{\kappa W_N(\xi,\sigma)}\nu_-(d\xi)\leq e^{4 D \kappa^2 \beta^2 C_1(\alpha)}. 
\end{equation}
Changing $\kappa\to -\kappa$ in \eqref{UB for num of f^N_+} we also obtain
\begin{equation}\label{minusUB for num of f^N_+}
    \int_{X_-} e^{-\kappa W_N(\xi,\sigma)}\nu_-(d\xi)\leq e^{4 D \kappa^2 \beta^2 C_1(\alpha)}. 
\end{equation}
Thus by applying the Cauchy-Schwarz inequality, we obtain
\begin{eqnarray*}
1 &=& \biggl[ \int_{X_-} e^{\kappa W_N(\xi,\sigma)/2}\,e^{-\kappa W_N(\xi,\sigma)/2}\,\nu_-(d\xi)\biggr]^2
\;\le\; \int_{X_-} e^{-\kappa W_N(\xi,\sigma)}\nu_-(d\xi)\, \int_{X_-} e^{\kappa W_N(\xi,\sigma)}\nu_-(d\xi)\\
&\le& e^{4 D \kappa^2 \beta^2 C_1(\alpha)} \cdot \int_{X_-} e^{\kappa W_N(\xi,\sigma)}\nu_-(d\xi)
\end{eqnarray*}
which yields the lower bound 
\begin{equation}\label{LB for num of f^N_+}
    e^{-4 D \kappa^2 \beta^2 C_1(\alpha)}\leq \int_{X_-} e^{\kappa W_N(\xi,\sigma)}\nu_-(d\xi) 
\end{equation}
for all $\sigma \in X_+$, $N\in\N$.  Putting \eqref{UB for num of f^N_+} and \eqref{LB for num of f^N_+} together yields the bounds
\begin{equation}\label{bounds for denom of f^N and f^N_+}
    e^{-4 D \kappa^2 \beta^2 C_1(\alpha)}
    \leq \int_{X_+}\int_{X_-} e^{\kappa W_N(\xi,\sigma)}\nu_-(d\xi)\nu_+(d\sigma)
    \leq e^{4 D \kappa^2 \beta^2 C_1(\alpha)} 
\end{equation}
which implies that the family $\{f^{[N]}_+:N\in\N\}$ is uniformly bounded from above and below:
\begin{equation}\label{uniform ub lb boundedness of half-line densities}
    e^{-8 D \beta^2 C_1(\alpha)}
    \leq f^{[N]}_+(\sigma)
    \leq e^{8 D \beta^2 C_1(\alpha)}\;.
\end{equation}

\noindent
\textit{Equicontinuity:}  
As the denominator $\int_X e^{-W_N}d\nu^{(0)}$ is uniformly bounded from above and below as shown in (\ref{bounds for denom of f^N and f^N_+}), it is enough to show that the family $\Big\{\int_{X_-}e^{-W_N(\xi,\cdot)}\nu_-(d\xi): N\in\N\Big\}$  is equicontinuous. 
Consider $n\in\N$ and configurations $\sigma, \tilde \sigma\in X_+$ such that $\sigma_0^{n-1}=\tilde\sigma _0^{n-1}$.
Then 
\begin{equation}
    \Big|\int_{X_-}\Big[e^{-W_N(\xi,\sigma)}-e^{-W_N(\xi,\tilde\sigma)} \Big]d\nu_-(d\xi)\Big|
    \leq \int_{X_-} e^{-W_N(\xi,\sigma)}\cdot\Big | e^{W_N(\xi,\sigma)-W_N(\xi,\tilde\sigma)}-1 \Big|\nu_-(d\xi).
\end{equation}
Thus, by the Cauchy-Schwarz inequality,
\begin{equation}
    RHS
    \leq \Big(\int_{X_-}e^{-2W_N(\xi,\sigma)}\nu_-(d\xi) \Big)^{\frac{1}{2}} \Big(\int_{X_-}\Big[ e^{W_N(\xi,\sigma)-W_N(\xi,\tilde\sigma)}-1 \Big]^2 \Big)^{\frac{1}{2}}
\end{equation}
and (\ref{UB for num of f^N_+}) yields that
\begin{equation}\label{equicontinuity important bound}
    \Big|\int_{X_-}\Big[e^{-W_N(\xi,\sigma)}-e^{-W_N(\xi,\tilde\sigma)} \Big]d\nu_-(d\xi)\Big|
    \leq C_2 \Big(\int_{X_-}\Big[ e^{W_N(\xi,\sigma)-W_N(\xi,\tilde\sigma)}-1 \Big]^2 \Big)^{\frac{1}{2}}
\end{equation}
with $C_2:= e^{8 D  \beta^2 C_1(\alpha)}$.  We will bound the last integral by bounding the exponent.
\smallskip

Note that
\begin{equation}
    W_N(\xi,\sigma)-W_N(\xi,\tilde\sigma)
    =-\beta\sum_{j=n}^N(\sigma_j-\tilde\sigma_j)\sum_{i=1}^{N}\frac{\xi_{-i}}{(i+j)^{\alpha}}
\end{equation}
and, thus, for all $k\in\N$,
\begin{equation}
    \delta_{-k}(W_N(\cdot,\sigma)-W_N(\cdot,\tilde\sigma))
    =2\beta\Big|\sum_{j=n}^N\frac{\sigma_j-\tilde\sigma_j}{(k+j)^{\alpha}} \Big|
    \leq 4\beta \sum_{j=n}^N \frac{1}{(k+j)^{\alpha}}\;.
\end{equation}
Hence, for sufficiently large $n$ and $N>n$,
\begin{equation}\label{oscillation estimate for the difference of W_N's}
    ||\underline\delta(W_N(\cdot,\sigma)-W_N(\cdot,\tilde\sigma))||_2^2
    \leq 16\beta^2\sum_{k=1}^{\infty}\Big(\sum_{j=n}^N\frac{1}{(k+j)^{\alpha}} \Big)^2
    \leq 32 \beta^2 \sum_{k=n+1}^{\infty}\frac{1}{k^{2(\alpha-1)}}=:u_n.
\end{equation}
with 
\[
\lim\limits_{n\to\infty}u_n=0 \;\mbox{ if and only if }\; \alpha>\frac{3}{2}\;.
\]
The spin-flip invariance of $\nu_-$ implies that for all $N\in\N$,
\[
\int_{X_-} [W_N(\xi,\sigma)-W_N(\xi,\tilde\sigma)]\,\nu_-(d\xi)=0\;.
\]
Then
since $\Phi^-$ satisfies (\ref{DUC for general interactions}), we have, from the second part of Theorem \ref{GCB and MCB}, that for all $m\in\N$,
\begin{eqnarray}\label{apply MCB to difference of W_N}
    \int_{X_-}\Big|W_N(\xi,\sigma)-W_N(\xi,\tilde\sigma) \Big|^m\nu_-(d\xi)
   & \leq& \Big(\frac{D  ||\underline\delta(W_N(\cdot,\sigma)-W_N(\cdot,\tilde\sigma))||_2^2   }{2} \Big)^{\frac{m}{2}}m\,\Gamma\Big (\frac{m}{2}\Big) \nonumber\\
   &=&  m\, v_n^m\, \Gamma \Big(\frac{m}{2}\Big)
\end{eqnarray}
with
\[
v_n:=\Big(\frac{D\,u_n}{2} \Big)^{\frac{1}{2}}\;.
\]
We conclude by expanding the square in the right-hand side in (\ref{equicontinuity important bound}):

Expanding the square in the right-hand side of \eqref{equicontinuity important bound},
\begin{eqnarray}
 \lefteqn{\int_{X_-}\Big[ e^{W_N(\xi,\sigma)-W_N(\xi,\tilde\sigma)}-1 \Big]^2 \nu_-(d\xi)}\nonumber\\
 &=& \int_{X_-}\Big[ e^{2\left[W_N(\xi,\sigma)-W_N(\xi,\tilde\sigma)\right]}- 2\, e^{W_N(\xi,\sigma)-W_N(\xi,\tilde\sigma)} 
  +1 \Big]\nu_-(d\xi)\nonumber\\
 &\le& 1+ \sum_{m=1}^{\infty} \frac{2^m+2}{m!} \int_{X_-}\Bigl|W_N(\xi,\sigma)-W_N(\xi,\tilde\sigma) \Bigr|^m\nu_-(d\xi).
\end{eqnarray}
Through elementary analysis one can show that
\begin{equation}
    \sum_{m=1}^{\infty}\frac{2^m+2}{m!} m \,v_n^m\, \Gamma \Big(\frac{m}{2}\Big)
    \leq 4v_n+3 e^{v_n^2}+2e^{4v_n^2}-5.
\end{equation}
Therefore, we obtain from (\ref{equicontinuity important bound}) and \eqref{apply MCB to difference of W_N}, 
\begin{equation}
    \Big|\int_{X_-}\Big[e^{-W_N(\xi,\sigma)}-e^{-W_N(\xi,\tilde\sigma)} \Big]d\nu_-(d\xi)\Big|
    \leq C_2 ( 4v_n+3 e^{v_n^2}+2e^{4v_n^2}-5)^{\frac{1}{2}}.
\end{equation}
Since $\lim\limits_{n\to\infty}v_n=0$, we conclude that the family $\Big\{\int_{X_-}e^{-W_N(\xi,\cdot)}\nu_-(d\xi): N\in\N\Big\}$  is indeed equicontinuous.
\end{proof}
\section{Final Remarks and Future Directions}\label{final remarks}
The main results of the present work rely on two major assumptions:
uniqueness of the Gibbs measure and validity of the Gaussian concentration inequalities for that measure.  We informally refer to the combination of these two conditions as \emph{strong uniqueness}. Strong uniqueness holds for a wide class of Gibbs interactions (potentials).
 
However, we would like to end with a discussion of one particular model --
the Dyson model, which served as the primary motivation for the present work. 
The picture below (Figure \ref{fig:diagram}) summarises the current state of the art in the eigenfunction problem for the Dyson potential with the parameters $\alpha>1$, $\beta>0$.

\begin{figure}[H]
\centering
\begin{center}
\begin{tikzpicture}
\begin{axis}[
    width=.9\textwidth,
    height=.43\textheight,
  axis lines = left,
  xlabel = $\alpha$,
  ylabel = $\beta$,
  xlabel style={
    at={(current axis.right of origin)},
    anchor=north west,
  },
  ylabel style={
    at={(current axis.above origin)},
    anchor=south east,
    rotate=270, 
  },
  xmin=1,
  xmax=3.0,
  ymin=0, 
  ymax=20, 
  xtick={1, 2}, 
  ytick distance=30, 
  extra x ticks={1.5,1.3},
  extra x tick labels={$3/2$, $\alpha'$},
  extra y ticks={15.5},
  extra y tick labels={$\beta'$},
]

\addplot [red, domain=1:2]{e^((1.48)*x)-e^(1.48)};
\addlegendentry{\tiny{$\beta_c(\alpha)$}}

\addplot [blue, domain=1.5:2]{x^(3.44) - 0.5};
\addlegendentry{\tiny{$\beta_*(\alpha)$}}

\addplot [blue, dashed, domain=1:2]{x^(3.3) -1};
\addlegendentry{\tiny{$\beta_c^{\tiny{\text{\textit{DU}}}}(\alpha)$}}

\node[] at (axis cs: 2.5,10.0) {\tiny{\textbf{(a) }}};

\addplot [white, name path=A,domain=2:2.96] {19.44};
\addplot [name path=B,domain=2:2.96] {0};
\addplot [green!20] fill between [of=A and B];

\node[] at (axis cs: 1.85,2.8) {\tiny{\textbf{(b)}}};

\addplot [blue, dashed, name path=A,domain=1.5:2] {x^(3.44) - 0.5};
\addplot [name path=B,domain=1.5:2] {0};
\addplot [green!20] fill between [of=A and B];

\node[] at (axis cs: 1.34,0.6) {\tiny{\textbf{(c)}}};

\addplot [blue, dashed, name path=A,domain=1:1.5] {x^(3.3) -1};
\addplot [name path=B,domain=1:1.5] {0};
\addplot [green!9] fill between [of=A and B];

\node[] at (axis cs: 1.42,3.0) {\tiny{\textbf{(d)}}};

\addplot [red,name path=A,domain=1:1.5] {e^((1.48)*x)-e^(1.48)};
\addplot [blue,dashed, name path=B,domain=1:1.5] {x^(3.3) -1};
\addplot [yellow!20] fill between [of=A and B];

\node[] at (axis cs: 1.4,10.3) {\tiny{\textbf{(f)}}};

\addplot [white, name path=A,domain=1:1.5] {19.44};
\addplot [red, name path=B,domain=1:1.5] {e^((1.48)*x)-e^(1.48)};
\addplot [red!15] fill between [of=A and B];

\addplot [red!17, name path=A,domain=1.5:2] {15.5};
\addplot [red, name path=B,domain=1.5:2] {e^((1.48)*x)-e^(1.48)};
\addplot [red!15] fill between [of=A and B];

\addplot [white, name path=A,domain=1.5:1.7] {19.44};
\addplot [red!17, name path=B,domain=1.5:1.7] {14};
\addplot [red!15] fill between [of=A and B];

\node[] at (axis cs: 1.67,17.3){\tiny{\textbf{(g)}}};

\addplot [white, name path=A,domain=1.3:2] {19.44};
\addplot [red!15, name path=B,domain=1.3:2] {15.5};
\addplot [red!11] fill between [of=A and B];

\node[] at (axis cs: 1.8,8.3) {\tiny{\textbf{(e)}}};

\addplot [red,name path=A,domain=1.5:2] {e^((1.48)*x)-e^(1.48)};
\addplot [blue, dashed, name path=B,domain=1.5:2] {x^(3.44) - 0.5};
\addplot [yellow!20] fill between [of=A and B];

\draw [dashed, black!100] (1.5,0) -- (1.5,19.44);
\draw [thick, dashed, black!100] (2,0) -- (2,19.44);
\addplot[mark=none, black, dotted] coordinates {(1.3,15.5) (1.3,19.44)};
\addplot[mark=none, black, dotted] coordinates {(1.3,15.5) (2,15.5)};

\end{axis}
\end{tikzpicture}
\end{center}

\caption{Eigenfunctions for the Dyson model across the phase diagram}
\label{fig:diagram}
\end{figure}
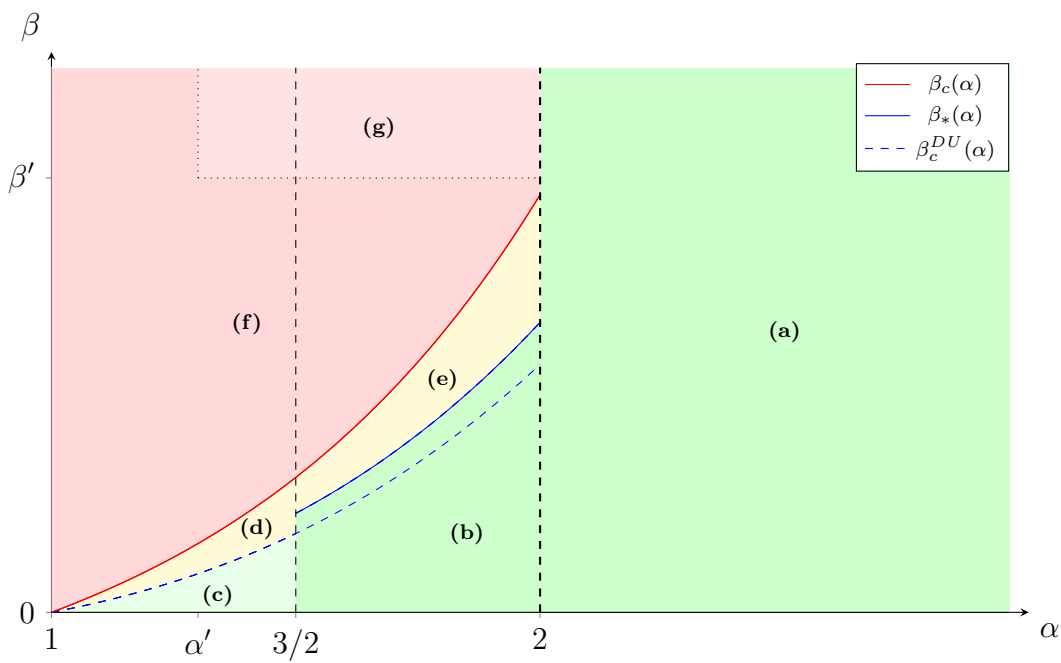
\begin{itemize}
    \item[(a)] 
For $\alpha>2$, at all temperatures the Dyson potential has summable variation, and hence
the classical results of Walters \cites{Walters1975,Walters1978} allow one to conclude that
 the transfer operator admits a continuous eigenfunction with summable variation.
 Note also that for $\alpha>2$ we also have strong uniqueness: uniqueness is due to Bowen
  \cites{Bow}, and the Gaussian Concentration Bounds have been established in 
  \cites{CMU2020,CGT2023}. Thus our results are also applicable, although summable variation 
  of the unique continuous positive eigenfunction requires a separate argument.

    \item[(b)] For $\alpha \in (3/2,2]$ existence of a continuous
    eigenfunction was first established by Johansson, \"Oberg, and Pollicott for sufficiently small $\beta< \beta_*(\alpha)$, see \cite{JOP2023} for details.
    Theorem C also applies in this region for $\beta<\beta_c^{DU}(\alpha)$, i.e., $\beta$'s satisfying Dobrushin's
    uniqueness condition. 
    Note $\beta_c^{DU}(\alpha)\le \beta_*(\alpha)$, thus the result of Theorem C might be weaker. However, the proof in  \cite{JOP2023}  relies on rather specific properties of the Dyson model, e.g., representation of the model via the random cluster model which typically restricts one to ferromagnetic models.
    
    \item[(c)] We support the conjecture by Johansson, \"Oberg, and Pollicott, that a positive continuous eigenfunction should exist for all $\beta$ below the critical value $\beta_c(\alpha)$.

    \item[(d)]  In this region, we have established the existence of an integrable eigenfunction.
     (c.f. the second part of Theorem \ref{main theorem about the Dyson model}). 
     We conjecture that the result is sharp: for $\alpha\le 3/2$, the transfer operator does not have a continuous eigenfunction.
    
    \item[(e)]  Similar to the case $\alpha>3/2$, we conjecture that the result must hold for all $\beta$'s below the critical value $\beta_c(\alpha)$.
    
    \item[(f)] In the supercritical regime: $\alpha>1$ and $\beta>\beta_c(\alpha)$, we believe that transfer operators do not admit integrable eigenfunctions.
    
    \item[(g)] The only result which applies to the supercritical phase is  \cite{BEvELN}, where 
    it has been shown that in a particular region, both pure phases of the Dyson model $\Phi$ are not $g-$measures. That immediately implies that
   the transfer operator does not have a continuous eigenfunction; otherwise, the normalized function $g=\frac{h\cdot e^{\phi}}{\lambda\cdot h\circ S}$ would become a \textit{$g-$function} for all the phases, in particular, the pure phases of the model.  
\end{itemize}
To summarize the picture, we conjecture that for the Dyson potential,
transfer operator does not have an eigenfunction in the supercritical regime,
and does have an eigenfunction in the subcritical regime. The smoothness
of the eigenfunction is varying with $\alpha$: from
summable variation for $\alpha>2$,  to continuous for $\alpha \in (3/2,2]$, to
$L^1$ but not continuous for $\alpha\le 3/2$.

The key to establishing  properties of transfer operators for the Dyson potential is, in our opinion,
a proper understanding of the probabilistic properties of the \emph{left-right interaction energy} function:
 $$
 W(\xi,\sigma)=\sum_{i<0\leq j}-\frac{\beta \xi_i\sigma_j}{(j-i)^\alpha},  \quad\xi\in X_-,\, \sigma\in X_+.
 $$
 For example, in the case $\alpha>2$, the results follow almost immediately from the simple observation that the interaction energy function is uniformly bounded
 $$
 \sup_{\xi\in X_-,\sigma\in X_+} |W(\xi,\sigma)|=\beta \sum_{i<0\leq j} \frac 1{(j-i)^\alpha}<\infty,
 $$
 and all expressions for densities above automatically lead to continuous functions.
 The next interesting 'critical value' $\alpha=3/2$  also appears quite naturally:
 The   left-right interaction energies $W(\xi,\sigma)$ for a fixed $\sigma\in X_+$, but $\nu_-$-random configurations $\xi\in X_-$, can be represented as
\[
W(\xi,\sigma)=\sum_{i<0}  Z_i , \quad Z_i=\xi_i\Bigl( -\sum_{j\ge 0} \frac{\sigma_j}{(j-i)^\alpha}\Bigr)
.
\]
The variances $\textsf{var}(Z_i)=\mathcal O(|i|^{-2(\alpha-1)})$  become summable for $\alpha>3/2$ [c.f. \eqref{eq:r-variance}].  
Hence, assuming weak correlations, the condition $\alpha>3/2$ corresponds to the almost sure existence of the left-right interaction energy for random conditions on the left half-line interacting with any (in particular, the all plus or all minus) configuration on the right half-line.
The concentration inequality can be interpreted as the rigorous transcription of this observation, and the Dobrushin condition as the guarantor of weak correlations.
 
 We strongly believe that the analysis of the left-right interaction energy function can be extended to
 the whole subcritical regime $\beta<\beta_c(\alpha)$.

\medskip

We finish the discussion with two interesting questions.
As customary in dynamical systems, we study continuous potentials $\phi\in C(X_+)$.
In Section \ref{From $C(X_+)$ to the space of interactions on $X$}, however, we switch to the language
of Statistical Mechanics and assume that the potential can be represented as
$\phi=-\sum_{0\in V\Subset\Z_+}\Phi_V$ for some translation invariant UAC interaction $\Phi$ on $Z$, c.f. \cite{Rue}. This is clearly the case for the Dyson potential. In general, however, we 
do not know any reasonable description of the class of such potentials $\phi$.

Finally, in the opposite direction, under which conditions, is a Gibbsian specification $\gamma$ on the half-lattice $\Z_+$ can be represented as $\gamma=\gamma^{\phi}$ for some $\phi\in C(X_+)$.
One possible approach to finding such representations would be extending  the Kozlov-Sullivan characterization on $\Z$ 
to the half-line $\Z_+$ \cite{BGMMT2020}.

\section*{Acknowledgements}
The authors would like to thank Anders \"Oberg for his useful correspondence, Frank Redig for helpful advice on concentration inequalities,  and A.v.E. thanks Eric Endo and Arnaud Le Ny for an earlier collaboration which provided some helpful background results.

\end{document}